\documentclass[preprint,12pt]{elsarticle}

\usepackage[T1]{fontenc}
\usepackage{amsmath,amssymb}
\usepackage{graphicx}
\usepackage{booktabs}
\usepackage{tabularx}
\usepackage{array}
\newcolumntype{Y}{>{\raggedright\arraybackslash}X}
\usepackage{listings}
\usepackage{xcolor}
\usepackage{url}
\usepackage{hyperref}
\newcommand{\odatse}{ODAT-SE}
\newcommand{\twodmat}{2DMAT}

\newcommand{\code}[1]{\texttt{#1}}
\newcommand{\angstrom}{\text{\AA}}

\journal{Computer Physics Communications}

\begin{document}

\begin{frontmatter}

\title{\odatse{}: A modular and parallel framework for simulation-based inverse problems}

\author[issp]{Tatsumi Aoyama\corref{cor1}}
\ead{aoym@issp.u-tokyo.ac.jp}
\author[issp]{Kazuyoshi Yoshimi\corref{cor1}}
\ead{k-yoshimi@issp.u-tokyo.ac.jp}
\author[issp]{Yuichi Motoyama}
\author[sokendai,nifs]{Akito Nakano}
\author[issp]{Hengyu Li}
\author[nifs]{Yusuke Konishi}
\author[nifs,sokendai]{Takeo Hoshi}
\cortext[cor1]{Corresponding authors}
\address[issp]{Institute for Solid State Physics, The University of Tokyo, Kashiwa, Chiba 277-8581, Japan}
\address[sokendai]{The Graduate University for Advanced Studies (SOKENDAI), Toki, Gifu 509-5292, Japan}
\address[nifs]{National Institute for Fusion Science, Toki, Gifu 509-5292, Japan}

\begin{abstract}
Open Data Analysis Tool for Science and Engineering (\odatse{}) is an
open-source Python framework for inverse problems that are solved by
repeated forward simulations.
It succeeds \twodmat{}, which was developed for surface-structure
analysis with quantum-beam diffraction and has since been applied to
other fields.
In version~3, the search algorithms were separated from the
application-specific solvers, which became independent packages, and
the software was renamed; this paper describes the current
version~4.0.
The framework connects optimization and sampling algorithms to
solvers through a common interface for parameter mapping, constraints,
and objective evaluation, so that an existing simulation code can be
used with any of the algorithms.
Relative to the previously published \twodmat{}, it provides
two-level MPI parallelism, which distributes candidate evaluations
over processes and each evaluation over a parallel solver, checkpoint
and continuation facilities, random search, global optimizers from
SciPy, tensor-train optimization, and a Bayesian workflow that
estimates partition-function ratios, model evidence, effective noise
levels, and marginal posterior distributions.
We validate the Bayesian workflow against closed-form results for a
linear regression, benchmark tensor-train optimization on continuous,
binary, and mixed-integer problems, and verify the two-level
parallelization by the same optimum for all process decompositions.
Taking up the positron-diffraction example of the \twodmat{} paper,
the workflow identifies the correlated structural solutions and an
effective noise scale close to the imposed perturbation, while
exposing incomplete convergence of symmetry-related mode weights.
Solver packages, templates, and a gallery of examples support the
use of the framework as a general platform for inverse-problem
analysis.

\medskip
\noindent \textbf{PROGRAM SUMMARY}

\begin{small}
\noindent
{\em Program Title:} Open Data Analysis Tool for Science and
Engineering (\odatse{}) \\
{\em CPC Library link to program files:} (to be added by Technical Editor) \\
{\em Developer's repository link:}
\url{https://github.com/issp-center-dev/ODAT-SE} \\
{\em Licensing provisions:} MPL-2.0 \\
{\em Programming language:} Python 3 (3.9--3.14) \\
{\em Supplementary material:} Manual, API reference, and
Gallery tutorials; see the Data availability statement. \\
{\em Journal reference of previous version:}*
Comput.\ Phys.\ Commun.\ 280 (2022) 108465 \\
{\em Does the new version supersede the previous version?:}* Yes. \\
{\em Reasons for the new version:}*
Applications beyond surface diffraction require reusable solver
interfaces, MPI-parallel forward evaluations, restartable execution,
discrete optimization, and integrated Bayesian analysis. \\
{\em Summary of revisions:}*
Renaming and modularization; independent solver packages; two-level
MPI parallelism; checkpoint and continuation; random, SciPy global,
and tensor-train optimization; Monte Carlo extensions; evidence and
posterior post-processing; license changed to MPL-2.0. \\
{\em Nature of problem:}
Inverse problems in which physical parameters are inferred from
measured data by comparing the data with the prediction of a forward
model, such as the determination of atomic positions at a surface from
quantum-beam diffraction intensities.
The forward model is often an existing simulation program whose single
evaluation is expensive, the objective function may have several
nearly degenerate minima and strongly correlated parameter directions,
and the noise level of the data is usually unknown.
Beyond a best-fit solution, the analysis therefore requires the
distribution of the parameters, the identification of competing
solutions, and an estimate of the effective noise level, obtained with
a budget of forward evaluations that can be distributed over many
processes and spread over restartable runs on batch systems. \\
{\em Solution method:}
The framework separates the search algorithm from the direct problem
solver through a common interface, so that any algorithm can be
combined with any solver, including external simulation programs and
objective functions written in Python.
The algorithms comprise local optimization (Nelder--Mead and other
SciPy minimizers with optional basin hopping), grid search, Bayesian
optimization, random and quasi-random search, global optimization
(differential evolution, SHGO, DIRECT, and dual annealing), tensor-train
optimization on discretized and mixed-variable spaces, and replica
exchange and population annealing Monte Carlo.
The Monte Carlo samplers estimate partition-function ratios along the
temperature schedule; post-processing tools evaluate the model
evidence as a function of the noise precision, select its maximum as
an empirical-Bayes estimate of the noise level, and construct marginal
posterior distributions from the samples at that temperature.
Candidate evaluations are distributed over MPI processes, and each
evaluation can itself be executed by a parallel solver; supported
algorithms can be checkpointed, resumed, and continued. \\
{\em Additional comments including restrictions and unusual features:}
External routines/libraries: NumPy and matplotlib; tomli for Python $<3.11$.
Feature-specific dependencies include SciPy, mpi4py, PHYSBO, and tqdm.
Solver packages may require external simulators.
Nonlinear constraints require solver-level treatment.
\code{minsearch} and \code{global\_search} do not support restart.
Solver-level parallelism requires a compatible parallel solver
implementation.

\medskip\noindent
* Items marked with an asterisk are only required for new versions
of programs previously published in the CPC Program Library.
\end{small}
\end{abstract}

\end{frontmatter}

\section{Introduction}
\label{sec:intro}

Inferring physical parameters from measured data is a central task
in experimental science. When a forward model predicts observables
$D_{\mathrm{cal}}(x)$ for parameters $x=(x_1,\ldots,x_N)$, parameter
estimation can be formulated as minimizing
\begin{equation}
  f(x)=d\bigl(D_{\mathrm{cal}}(x),D\bigr),
  \label{eq:objective}
\end{equation}
where $D$ denotes the reference data and $d$ measures their discrepancy
from the prediction. We call the program evaluating the prediction
the \emph{direct problem solver}, following \twodmat{} terminology.
Beyond a single optimum, physical interpretation often requires
parameter uncertainties, correlations, and identification of competing
solutions. This combination of simulation and statistical inference
is part of measurement informatics~\cite{Nagata2012}.

The main computational expense frequently lies in the forward
simulation rather than the search algorithm. A practical framework
must therefore reuse existing physical codes, distribute repeated
evaluations efficiently, and support long calculations on batch systems.
These requirements become particularly important when the solver is
an external program, or itself uses MPI.
Quantum-beam surface diffraction provides a representative example.
TRHEPD~\cite{Hugenschmidt2016,Fukaya2019},
SXRD~\cite{Feidenhansl1989,Tajiri2020}, and
LEED~\cite{VanHove1993} determine near-surface atomic coordinates
by comparing intensities with diffraction simulations.
Their objective landscapes can contain nearly degenerate minima
and correlated parameter directions~\cite{Hoshi2022}.
Extended-ensemble Monte Carlo methods help characterize these
ambiguities beyond local fitting~\cite{Anada2017,Nagai2020}.

To facilitate such analyses, some of the authors developed \twodmat{}~\cite{Motoyama2022,Tanaka2020,Tanaka2023}, combining
diffraction solvers with local optimization, grid search, Bayesian
optimization, replica exchange, and population annealing.
Applications subsequently extended to surface phases~\cite{Tsujikawa2022,Sato2025},
magnetotransport~\cite{Liu2025}, fusion-plasma
diagnostics~\cite{Sakai2026jinst,Sakai2026ppcf,Kobayashi2026,Kin2026pfr},
magnetic-configuration searches~\cite{Kin2026ppcf}, and performance
modeling~\cite{Kohashi2022}.
Related developments include a fast diffraction simulator designed
for global search~\cite{Kudo2024} and the co-evolution of
PHYSBO~\cite{Motoyama2022physbo,Motoyama2026physbo}; the software is also
distributed through MateriAppsLIVE!~\cite{Motoyama2022materiapps}.
These applications motivated a
domain-independent core and independently maintained solver interfaces.
Accordingly, the framework was reorganized in version~3.0, released in
2024: the search algorithms were separated from the application-specific
solvers, which became independent packages, and the software was renamed
Open Data Analysis Tool for Science and Engineering (\odatse{}).
The version numbering continues from \twodmat{}, so that the first
release of \odatse{} is v3.0.
At the same time, the license was changed from the GNU General Public
License version~3 (GPL-3.0) to the Mozilla Public License version~2.0
(MPL-2.0).
The GPL imposes strong constraints when the framework is embedded as a
library in other software, whereas MPL-2.0 retains a file-level copyleft
and is better suited to use as a library.

Existing packages address related needs: emcee~\cite{ForemanMackey2013}
and PyMC~\cite{AbrilPla2023} provide posterior sampling,
Optuna~\cite{Akiba2019} supports black-box optimization, and
refnx~\cite{NelsonPrescott2019} and GenX~\cite{GlavicBjorck2022}
combine reflectometry models with fitting and uncertainty analysis.
The focus of \odatse{} is the integration of interchangeable search
and sampling methods with application-specific simulation codes,
including compatible MPI-parallel solvers, through a common
configuration and execution interface.

This paper presents \odatse{} on the basis of its current version~4.0
and describes the developments since the previously published
\twodmat{} v2.
The main developments are a modular solver architecture; two-level
parallelism with restart support for selected algorithms; additional
continuous, discrete, and mixed-variable optimization methods; and
a workflow connecting Monte Carlo normalization estimates to
evidence-based noise estimation and posterior analysis.
The contribution is their integration for simulation-based
inverse analysis, rather than a new derivation of the underlying
optimization or sampling methods.

Section~\ref{sec:overview} gives an overview of \odatse{} with its
installation and a quick guide.
Section~\ref{sec:features} summarizes the changes since \twodmat{} and
describes the architecture and methods,
Section~\ref{sec:applications} evaluates them using controlled
benchmarks and synthetic diffraction data, and
Section~\ref{sec:summary} concludes the paper.
Supporting examples and normalization details are collected in
the appendices.

\section{Overview of \odatse{}}
\label{sec:overview}

This section introduces \odatse{} from the user's point of view: what the
framework provides, how it is installed, and how a calculation is run.

\subsection{What \odatse{} provides}
\label{sec:whatis}
\odatse{} is a Python framework that applies a search algorithm to a
direct problem solver in order to find the parameters that minimize the
objective function of Eq.~(\ref{eq:objective}), or to sample their
distribution.
The algorithm and the solver are selected independently of each other.
The available algorithms are local optimization by the Nelder--Mead and
related methods, grid search, Bayesian optimization, replica exchange and
population annealing Monte Carlo, random search, global optimization, and
tensor-train optimization; they are listed in Table~\ref{tab:algorithms}.
The core package contains two solvers: a set of analytical test functions
and a wrapper for a user-supplied Python function.
Solvers for quantum-beam diffraction and spectroscopy experiments are
distributed as separate packages, and users can add solvers for their own
simulation codes (Sections~\ref{sec:framework} and~\ref{sec:ecosystem}).
Between the algorithm and the solver, an adapter layer handles parameter
mapping, constraints, and logging (Section~\ref{sec:framework}).

\odatse{} can be used in two ways.
The \code{odatse} command runs a calculation specified by an input file,
as described in Section~\ref{sec:quickguide}.
The same components are available as a Python library, so that a
calculation can also be driven from a user's script.
A manual with tutorials, an API reference, and a Gallery of worked
examples and solver templates accompany the software; their locations are
given in the Data availability statement.

\subsection{Installation}
\label{sec:install}
\odatse{} requires Python~3.9 or later with NumPy and matplotlib, and
tomli for Python versions earlier than~3.11.
The other dependencies are optional and are needed only for particular
features: SciPy~\cite{Virtanen2020} for \code{minsearch},
\code{global\_search}, and the quasi-random designs of
\code{random\_search}; mpi4py~\cite{Dalcin2021} for MPI-parallel execution;
PHYSBO~\cite{Motoyama2022physbo,Motoyama2026physbo} for \code{bayes};
and tqdm for progress display.
The package is distributed through PyPI and is installed with
\begin{lstlisting}[language={}]
$ python3 -m pip install ODAT-SE
\end{lstlisting}
Specifying \code{'ODAT-SE[all]'} instead installs the optional
dependencies as well.
The solver packages of Section~\ref{sec:ecosystem} and the external
simulators that they call are installed separately.
The source repository contains unit tests based on pytest and CTest
regression tests that compare algorithm, restart, and MPI outputs with
reference data; both are run in continuous integration.

\subsection{Quick guide}
\label{sec:quickguide}
A calculation is specified by an input file in TOML format.
The \code{[base]} section gives settings shared by the whole calculation,
such as the number of parameters and the output directory;
\code{[solver]} selects the solver and sets its parameters; and
\code{[algorithm]} selects the search algorithm and defines the search
space and the method-specific settings.
An optional \code{[runner]} section configures the adapter layer.
Listing~\ref{lst:input} is a complete input file that minimizes the
Himmelblau function~\cite{Himmelblau1972},
\begin{equation}
  f(x_1,x_2)=(x_1^2+x_2-11)^2+(x_1+x_2^2-7)^2,
  \label{eq:himmelblau}
\end{equation}
with the Nelder--Mead method~\cite{NelderMead1965}, starting from the origin in the box
$[-6,6]^2$.
This function has four minima with $f=0$, one of which is located at
$(3,2)$.

\begin{lstlisting}[language={},caption={Minimal input for the Himmelblau function with Nelder--Mead.},label={lst:input}]
[base]
dimension = 2
output_dir = "output"

[solver]
name = "analytical"
function_name = "himmelblau"

[algorithm]
name = "minsearch"
seed = 12345

[algorithm.param]
min_list = [-6.0, -6.0]
max_list = [ 6.0,  6.0]
initial_list = [0.0, 0.0]
\end{lstlisting}

The calculation is started with
\begin{lstlisting}[language={}]
$ odatse input.toml
\end{lstlisting}
Algorithms that support MPI parallelization are launched in the same way
through the MPI launcher, for example
\code{mpiexec -np 8 odatse input.toml}.
During the run, the parameters and the objective value of each
evaluation are written to the standard output.
The result is written to \code{res.txt} in the output directory:
\begin{lstlisting}[language={}]
fx = 4.2278370361994904e-08
x1 = 2.9999669562950175
x2 = 1.9999973389336225
\end{lstlisting}
The search has converged to the minimum at $(3,2)$ after 79 evaluations
of the objective function.
Files produced by each process are placed in subdirectories named after
the process rank; here, \code{output/0/SimplexData.txt} records the path
of the search.
Replacing the \code{[algorithm]} section selects another method with its
own search-space and method-specific settings, while the \code{[solver]}
section is unchanged.
Figure~\ref{fig:himmelblau} shows random and global searches of the
same function, and the manual provides a tutorial for each algorithm.

\section{New features}
\label{sec:features}
This section describes the developments since the previously published
\twodmat{} v2~\cite{Motoyama2022}.
Version~3.0 (November 2024) reorganized the code into the modular
framework of \code{Algorithm}, \code{Runner}, and \code{Solver}, moved
the experiment-specific solvers to separate packages, and introduced the
common checkpoint and restart interface.
Version~3.1 (April 2025) added the post-processing tools for population
annealing, and version~3.2 (November 2025) added random search and the
evaluation of partition functions in both Monte Carlo samplers.
Version~4.0 (August 2026) adds two-level MPI parallelization, global
optimization based on SciPy, tensor-train optimization, a choice of
local optimization methods, collective error handling, and the
installation of the post-processing tools as commands; the manual and
the tutorials were also revised and extended.
Version~4.0.1 (October 2026) is a maintenance release.
Table~\ref{tab:changes} summarizes these changes together with the
release that introduced each of them and the subsection that describes
it.

\begin{table}[tb]
\centering
\footnotesize
\caption{Changes from the previously published \twodmat{} v2 to \odatse{}.
Version numbers in parentheses give the release that introduced each change;
the last column gives the section describing it.}
\label{tab:changes}
\begin{tabularx}{\textwidth}{@{}lYYl@{}}
\toprule
Item & \twodmat{} v2 & \odatse{} & Sec. \\
\midrule
Parameter processing & Separate mapping and logging configuration & Mapping, constraints, and logging in \code{Runner} (v3.0) & \ref{sec:framework} \\
Parallelism & Algorithm-level MPI & Algorithm- and solver-level MPI (v4.0) & \ref{sec:mpi} \\
Restart & No common checkpoint interface & Resume and continuation for supported algorithms (v3.0) & \ref{sec:restart} \\
Algorithms & Five modules & Eight modules, adding random search (v3.2) and global and tensor-train search (v4.0) & \ref{sec:newalg} \\
Bayesian analysis & Posterior histograms & Partition-function ratios, evidence, noise estimation, and posterior tools (v3.1; extended in v3.2 and v4.0) & \ref{sec:bayes} \\
Distribution & \code{py2dmat}; bundled diffraction solvers & \code{ODAT-SE}; independent solver packages (v3.0) & \ref{sec:ecosystem} \\
License & GPL-3.0 & MPL-2.0 (v3.0) & -- \\
\bottomrule
\end{tabularx}
\end{table}

\subsection{Modular framework}
\label{sec:framework}
\begin{figure}[tb]
\centering
\includegraphics[width=\textwidth]{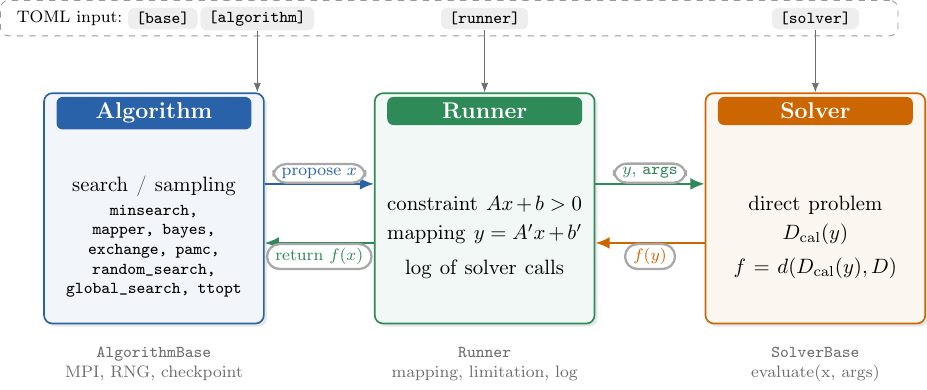}
\caption{Modular architecture of \odatse{}.
\code{Algorithm} proposes a search parameter vector $x$.
\code{Runner} checks the constraints, applies an optional mapping
to the solver parameter vector $y$, and requests an evaluation
from \code{Solver}.
The algorithm receives the returned objective value.
The TOML sections configure the corresponding layers, while
\code{[base]} provides global settings shared by the framework.}
\label{fig:layers}
\end{figure}

The search algorithms are independent of the physical problem, whereas
the forward simulation is specific to each problem and is often an
existing external program.
Applying the same inverse-analysis workflow to different problems
therefore requires that the two be developed and maintained
independently.
To this end, \odatse{} divides the evaluation of the objective function
among three components with distinct responsibilities,
\code{Algorithm}, \code{Runner}, and \code{Solver}
(Figure~\ref{fig:layers}); \code{Runner} is the adapter layer
mentioned in Section~\ref{sec:whatis}.

\code{Algorithm} explores the search coordinates $x$ using an optimization or sampling method.
\code{Runner} handles parameter constraints and optional mappings from the search coordinates $x$ to the solver coordinates $y$, and records evaluation logs.
\code{Solver} performs the application-specific calculation and returns a scalar objective value.
This separation allows different algorithms to use the same solver interface and keeps parameter processing outside the search and simulation implementations.

Algorithms derive from \code{AlgorithmBase}, which supplies common
MPI, random-number, timing, output, and checkpoint facilities.
They request evaluations through \code{self.runner.submit(x, args)},
where \code{args} contains auxiliary indices such as the step or
replica number.

The runner receives each evaluation point from the algorithm, checks
the constraints, applies the mapping, and calls the solver's
\code{evaluate} method.
Points violating the constraints are rejected by returning $+\infty$.
The built-in constraint is the componentwise linear inequality
$Ax+b>0$; arbitrary constraints can be imposed by a user-defined
\code{Limitation} class when the calculation is driven from a Python
script.
For admissible points, an optional mapping from the search coordinates
$x$ to the solver coordinates $y$ allows the search to use independent
or symmetry-adapted coordinates while the solver receives physical
parameters.
The built-in mapping is the affine mapping $y=A'x+b'$; arbitrary
mappings can likewise be supplied by a user-defined \code{Mapping}
class.
The dimensions of $x$ and $y$ may differ.

Application-specific solvers subclass \code{SolverBase} and implement
\code{evaluate(x, args)}. The method can evaluate a Python model
or prepare inputs, invoke an external simulator, and parse its output.
A driver constructs the solver, runner, and selected algorithm after
calling \code{odatse.initialize()}.
The built-in \code{function} solver wraps a Python callable, and the
\code{analytical} solver supplies test functions.
The same interface accommodates a fixed surrogate model
(\ref{app:surrogate}).
Detailed extension examples are provided by the solver templates
in Section~\ref{sec:ecosystem}.

\subsection{Two-level MPI parallelization}
\label{sec:mpi}

Simulation-based inverse analysis offers two opportunities for parallelization: evaluating multiple candidate parameter vectors concurrently and distributing the forward calculation for each candidate across multiple processes.
When individual forward calculations are computationally expensive, combining these two forms of parallelism allows computational resources to be allocated both across and within objective evaluations.

Building on the separation between algorithms and solvers described in Section~\ref{sec:framework}, \odatse{} supports two-level MPI parallelization on top of mpi4py~\cite{Dalcin2021}.
Algorithm-level parallelism coordinates concurrent candidate evaluations, while solver-level parallelism distributes each evaluation across a group of processes, provided that the solver is written to distribute its work over the group, as described below.

\begin{figure}[tb]
\centering
\includegraphics[width=\textwidth]{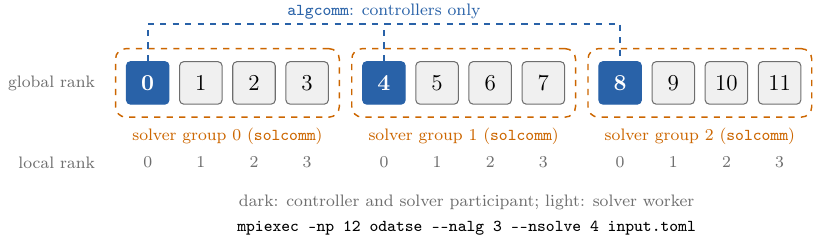}
\caption{Two-level MPI decomposition for
$n_{\mathrm{alg}}=3$ and $n_{\mathrm{solve}}=4$.
Each solver group has a separate communicator,
\code{solcomm}, containing one controller and three workers.
The controllers, shown as dark boxes at global ranks
0, 4, and 8, also form the algorithm communicator,
\code{algcomm}.
Numbers inside the boxes indicate global ranks; numbers
below the boxes indicate ranks within each solver group.
All processes in a group participate in the evaluation
of the same candidate parameter vector.}
\label{fig:comm}
\end{figure}

To implement this arrangement, the global communicator is split into $n_{\mathrm{alg}}$ solver groups, each containing $n_{\mathrm{solve}}$ processes, such that $N_{\mathrm{proc}} = n_{\mathrm{alg}} n_{\mathrm{solve}}$.
For example,
\begin{lstlisting}[language={}]
$ mpiexec -np 12 odatse --nalg 3 --nsolve 4 input.toml
\end{lstlisting}
creates three four-process groups (Figure~\ref{fig:comm}).
Within each group, the process of local rank 0 is the controller,
and the other members are workers that take part only in the solver
evaluations.
The controllers form the algorithm communicator, \code{algcomm}, on
which the search algorithm operates, and all members of a group form
its solver communicator, \code{solcomm}.
If one of the two options is omitted, it is inferred from the process
count, which must be divisible accordingly; without either option,
$n_{\mathrm{solve}}=1$.

After the runner checks constraints and applies mappings, the
controller broadcasts the solver parameters and auxiliary arguments
over \code{solcomm}. All group members call \code{evaluate}
with identical arguments.
Only the controller's return value is passed to the algorithm, and
those of the workers are discarded; the solver must therefore reduce
any partial results to local rank 0.

Solver-level parallelism is the responsibility of the solver: its
author distributes the work explicitly using
\code{odatse.mpi.solcomm()}, and the split does not by itself
parallelize a serial or external program.
Collective operations must be issued within \code{solcomm} and
called consistently by all members of the group; collectives over
the parent \code{MPI\_COMM\_WORLD} would block, because the
controllers and the workers are not synchronized outside the
evaluations.
Two further points are left to the solver author.
When thread parallelism is used in addition, the number of threads
must be matched to the cores assigned to each process.
When each solver process writes files, their names or directories
must be chosen so that the processes do not overwrite each other's
files.

The random-number generators are managed per algorithm rank so that
the parallel run is reproducible.
Each algorithm controller is seeded as
$s_r=s+r\Delta s$, where $s$ is \code{seed}, $r$ is its algorithm
rank, and $\Delta s$ (\code{seed\_delta}) defaults to 314159.
Keeping $n_{\mathrm{alg}}$ fixed therefore preserves the initial
algorithm random streams when $n_{\mathrm{solve}}$ changes;
bitwise reproducibility additionally depends on the solver and
floating-point reduction order.

The appropriate split balances available candidate concurrency
against solver scaling, as illustrated in Section~\ref{sec:app:parallel}.

\subsection{Checkpoint, restart, and error handling}
\label{sec:restart}

Sampling and optimization runs with expensive solvers can take hours
or days and are often executed on batch systems with limited job
times.
\odatse{} therefore provides a checkpoint mechanism that saves the
state of the algorithm periodically, so that an interrupted run can be
resumed and a completed run can be extended without repeating the
evaluations already made.
With \code{checkpoint = true}, each algorithm rank writes its state,
such as the current points, weights, and counters, together with its
random-number generator to a file (\code{status.pickle} by default)
at points where the state is consistent, for example at the end of a
Monte Carlo or annealing step.
The interval is given in steps (\code{checkpoint\_steps}) or in
seconds (\code{checkpoint\_interval}), and a checkpoint is always
written at the end of the run.
An external simulator called by the solver is not checkpointed: after
a restart, only the evaluations after the last checkpoint are
repeated.

The run mode is selected by a command-line option of \code{odatse},
parsed by \code{odatse.initialize()}.
The default \code{-{}-init} starts a new calculation.
\code{-{}-resume} restores an interrupted calculation, which must use
the same number of algorithm ranks, and \code{-{}-cont} extends a
completed run with modified algorithm parameters, for example
additional Monte Carlo steps, an extended annealing schedule, or a
larger optimization budget.
Table~\ref{tab:algorithms} lists which algorithms support resume and
continuation.
The framework checks the algorithm settings against the checkpoint
but not the solver: the objective, reference data, and solver
configuration must be kept unchanged by the user.
By default the random-number stream is restored; \code{-{}-reset\_rand}
reinitializes it instead, for example to make a continuation
statistically independent of the previous run.

Error handling in a parallel run is a related concern, because an
exception on one rank must not leave the other ranks waiting.
The execution of an algorithm consists of three stages, prepare, run,
and post.
At the end of each stage, the algorithm ranks exchange their status,
so that an exception on any rank, including a solver worker
inside \code{evaluate}, terminates the whole job in a
coordinated way with a nonzero exit status.
This is not general MPI fault tolerance: a failure inside a collective
operation cannot be recovered.
In addition, the runner can be instructed with
\code{ignore\_error} to treat a \code{RuntimeError} raised by the
solver as a failed evaluation, returning \code{NaN} for that point
and continuing the run.

\subsection{New and extended search algorithms}
\label{sec:newalg}

Table~\ref{tab:algorithms} lists the eight algorithm modules of
\odatse{}.
Five of them originate from \twodmat{}~\cite{Motoyama2022}; three are
new: \code{random\_search}, \code{global\_search}, and \code{ttopt}.
Among the original modules, \code{minsearch} and the two Monte Carlo
samplers have been extended, whereas \code{mapper} and \code{bayes},
which uses PHYSBO~\cite{Motoyama2022physbo,Motoyama2026physbo}, are
functionally unchanged apart from the restart support of
Section~\ref{sec:restart}.
The modules fall into three classes: optimization methods seek low
objective values, the Monte Carlo methods sample the
temperature-dependent distributions of Section~\ref{sec:bayes}, and
random search evaluates a set of points chosen without feedback from
the objective, for exploratory screening or as a baseline for the
other methods.
The new modules are described first, followed by the extensions.

\begin{table}[tb]
\centering
\footnotesize
\caption{Algorithm modules in \odatse{} v4.0.
Modules marked with $\dagger$ are new since \twodmat{}.
``MPI'' denotes parallel distribution of objective evaluations
over algorithm ranks, not solver-level parallelism.
R and C denote resume and continuation, respectively;
method-specific restrictions on continuation are given in the text.
All MPI execution requires mpi4py; additional feature-specific
dependencies are described in the text.}
\label{tab:algorithms}
\begin{tabularx}{\textwidth}{@{}lYcc@{}}
\toprule
Module & Methods or purpose & MPI & Restart \\
\midrule
\code{minsearch} & Nelder--Mead and other SciPy local minimizers; optional basin hopping & -- & -- \\
\code{mapper} & Exhaustive evaluation of grid or mesh points & yes & R \\
\code{bayes} & Bayesian optimization over candidate points using PHYSBO & yes & R, C \\
\code{exchange} & Replica exchange Monte Carlo & yes & R, C \\
\code{pamc} & Population annealing Monte Carlo & yes & R, C \\
\code{random\_search}$^\dagger$ & Pseudo-random, Sobol, Halton, and Latin hypercube sampling & yes & R, C \\
\code{global\_search}$^\dagger$ & SciPy global optimizers: differential evolution, SHGO, DIRECT, and dual annealing & DE, SHGO & -- \\
\code{ttopt}$^\dagger$ & Tensor-train optimization on discretized search spaces & yes & R, C \\
\bottomrule
\end{tabularx}
\end{table}

\paragraph{Random and global search}
\code{random\_search} evaluates a fixed number of points using
pseudo-random sampling, Sobol~\cite{Sobol1967} or
Halton~\cite{Halton1960} sequences, or Latin hypercube
sampling~\cite{McKay1979}, with quasi-random designs supplied by
SciPy~\cite{Virtanen2020}.
The scrambling of the sequences~\cite{Owen1998} is controlled by the
seed, so that the same point set is generated for any number of
algorithm ranks and a run is reproducible independently of the
process decomposition.
Sobol, Halton, and pseudo-random runs support continuation, because
further points of the same sequence can be appended;
Latin hypercube continuation is rejected because a design of a larger
size does not contain the smaller one.

\begin{figure}[tb]
\centering
\includegraphics[width=\textwidth]{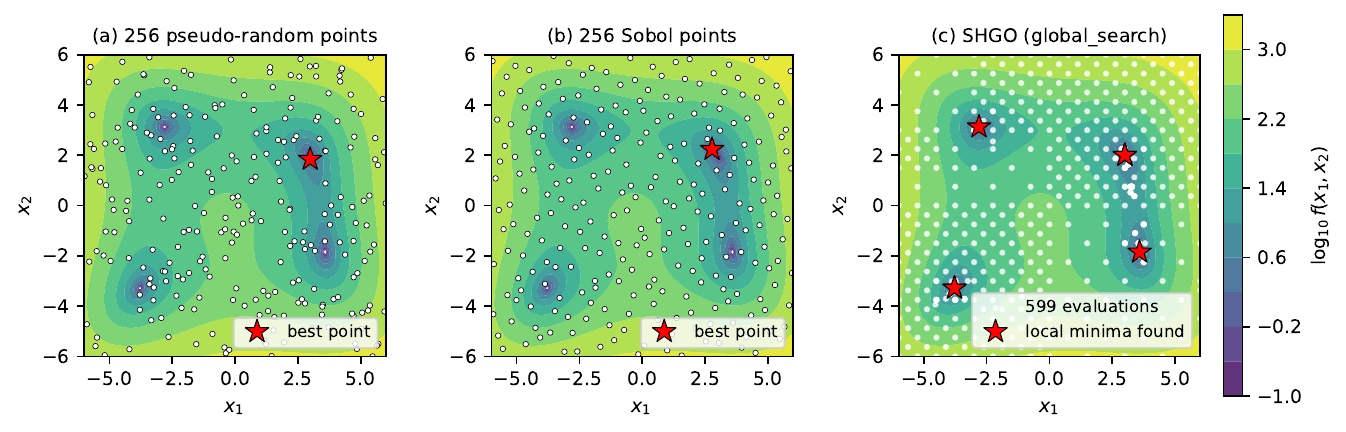}
\caption{Search of the Himmelblau function, Eq.~\eqref{eq:himmelblau},
on $[-6,6]^2$.
The background color represents $\log_{10}f$.
(a) Pseudo-random search with 256 points.
(b) A scrambled Sobol design with 256 points.
White dots indicate evaluated points, and the star in
each panel marks the smallest sampled objective value.
(c) SHGO with 128 sampling points per iteration and three iterations.
White dots indicate the 599 evaluated points, and stars
mark the four local minima found with $f<10^{-14}$.}
\label{fig:himmelblau}
\end{figure}

\code{global\_search} provides the global optimizers of SciPy:
differential evolution~\cite{StornPrice1997}, SHGO~\cite{Endres2018},
DIRECT~\cite{Jones1993}, and dual annealing~\cite{Xiang1997}.
The population of a differential-evolution generation and the
sampling points of an SHGO iteration can be evaluated at once and are
distributed over the algorithm ranks, whereas the local refinement of
SHGO and the iterations of DIRECT and dual annealing proceed sequentially and run on
algorithm rank 0, so that only solver-level parallelism applies to
them.
SHGO also records all local minima it discovers, which helps to
identify competing solutions of an inverse problem.
Figure~\ref{fig:himmelblau} illustrates these modules on the
Himmelblau function of Section~\ref{sec:quickguide}:
the scrambled Sobol design covers the search box more uniformly than
pseudo-random points with the same number of evaluations, and SHGO
locates all four minima with 599 evaluations.
The figure illustrates the interfaces, not a general performance
ranking.

\paragraph{Tensor-train optimization}
TTOpt~\cite{Sozykin2022} treats objective values on a grid as a tensor
without constructing the full tensor.
It uses tensor-train representations~\cite{Oseledets2011} and
cross approximation~\cite{OseledetsTyrtyshnikov2010} with
maximum-volume selection~\cite{Goreinov2010} to explore entries
and retain the smallest value encountered.
Axis $i$ has $N_i=P_i^{q_i}$ points, represented by $q_i$
indices of size $P_i$ through \code{p\_points} and
\code{q\_points}. Thus, a grid with $2^{20}$ points on each of
$N$ axes becomes a tensor with $20N$ binary indices; this
quantization keeps the tensor-train ranks small.

Controls include the maximum rank \code{r\_max}, the
requested-evaluation budget \code{max\_f\_eval}, and the tolerance of
the maximum-volume iteration \code{maxvol\_tol}.
Evaluation batches are distributed over algorithm ranks.
Repeated points are served from a cache, so that the number of
requested evaluations, which the budget counts, exceeds the number
of distinct solver calls; both are reported (Section~\ref{sec:app:ttopt}).
The sweep state and cache are checkpointed, and continuation
can increase the budget.
Continuous parameters are discretized, and a binary variable is
represented by a single index with $P_i=2$ and $q_i=1$.
This enables mixed-variable searches but does not guarantee
global optimality or automatically enforce nonlinear constraints,
which can be handled by penalty terms in the solver
(\ref{app:discrete}).

\paragraph{Local optimization}
\code{minsearch} extends Nelder--Mead~\cite{NelderMead1965} to
other SciPy methods selected by \code{method}, both derivative-free
methods and gradient-based methods with finite-difference gradients,
including Powell, COBYLA, COBYQA, L-BFGS-B, and SLSQP, subject to the
installed SciPy version.
For these methods, the search box is passed to SciPy as bounds,
whereas Nelder--Mead keeps the original behavior of returning
$+\infty$ outside the box.
An optional basin-hopping wrapper~\cite{WalesDoye1997} alternates
perturbations and local minimization.
Numerical gradients typically require approximately $N$
additional solver calls for forward differences or $2N$ for
central differences, which can be costly for physical simulations.
The internal states of the SciPy optimizers are not checkpointed,
hence the lack of restart support in Table~\ref{tab:algorithms}.

\paragraph{Monte Carlo extensions}
\label{sec:mcext}
Replica exchange~\cite{SwendsenWang1986,HukushimaNemoto1996}
and population annealing~\cite{HukushimaIba2003,Machta2010}
support continuous spaces with periodic axes and discrete meshes,
on which moves are proposed to neighboring points defined by a
neighbor list.
Both estimate partition-function ratios along the temperature
schedule.
For adjacent inverse temperatures,
\begin{equation}
  \frac{Z_{\beta_{k+1}}}{Z_{\beta_k}}
  =
  \left\langle
  e^{-(\beta_{k+1}-\beta_k)f(x)}
  \right\rangle_{\beta_k},
  \label{eq:zratio}
\end{equation}
where $\langle\cdot\rangle_{\beta_k}$ denotes the average over the
samples at $\beta_k$.
Replica exchange uses temperature-resolved sample averages;
population annealing uses annealing weights~\cite{Neal2001}.
The accumulated ratio is referenced to the initial inverse
temperature; population annealing can instead start from $\beta=0$
(\code{anneal\_from\_beta0}), which references the ratio to the
prior, as used in Section~\ref{sec:bayes}.

Population annealing additionally reports the participation ratio,
an effective number of walkers that indicates the concentration of
the weights.
For massively parallel runs, the per-rank output files can be
merged into combined files (\code{export\_combined\_files}), which
coarsens the granularity of the file I/O and reduces the load on the
file system; post-processing tools then expand the combined files
and collect the samples by temperature
(\code{odatse\_extract\_combined}, \code{odatse\_separateT}).
Diagnostics of population diversity are useful but do not alone
establish convergence between separated modes.

\subsection{Bayesian inference workflow}
\label{sec:bayes}

The Monte Carlo samplers of Section~\ref{sec:mcext} produce samples
and partition-function ratios along a temperature schedule.
This subsection describes how these outputs are turned into an
estimate of the noise level and into posterior distributions of the
parameters, which address the uncertainties and competing solutions
mentioned in Section~\ref{sec:intro}.
We first fix the statistical model.
For $K$ data subsets $D_\mu$ with $n_\mu$ points each, we assume
independent Gaussian errors with variance
$\sigma_\mu^2=1/(2\beta w_\mu)$.
The likelihood is then
\begin{equation}
  \pi(D\mid x;\beta)=
  \left[\prod_\mu
  \left(\frac{\beta w_\mu}{\pi}\right)^{n_\mu/2}\right]
  e^{-\beta f(x)},
  \label{eq:likelihood}
\end{equation}
where
\begin{equation}
  f(x)=\sum_{\mu=1}^{K}
  w_\mu\|D_{\mathrm{cal},\mu}(x)-D_\mu\|^2
  \label{eq:weighted_objective}
\end{equation}
is the objective function, $\|\cdot\|^2$ denotes the sum of squared
residuals over the $n_\mu$ points of a subset, and the weights
$w_\mu>0$ are normalized as $\sum_\mu w_\mu=1$, which fixes the
scale of $\beta$; for a single data set, $w_1=1$.
The inverse temperature $\beta$ of the samplers thus plays the role
of the noise precision.
The framework does not know the partition into subsets: the solver
must return $f(x)$ in the form of Eq.~\eqref{eq:weighted_objective}.

For a uniform prior on a finite admissible region $\Omega$,
\begin{equation}
  Z_\beta=\int_\Omega e^{-\beta f(x)}\,\mathrm{d}x,
  \qquad
  \pi(x\mid D;\beta)=\frac{e^{-\beta f(x)}}{Z_\beta}.
  \label{eq:posterior}
\end{equation}
Here $\Omega$ is the region actually sampled, that is, the search box
restricted by the constraints, and the prior is uniform in the search
coordinates $x$; if a mapping to solver coordinates is used, the
prior in the solver coordinates is uniform only if the mapping
preserves volume.
With $Z_0=V_\Omega$, the evidence is
\begin{equation}
  \log P(D;\beta)=
  \log\frac{Z_\beta}{Z_0}
  +\sum_\mu\frac{n_\mu}{2}
  \log\frac{\beta w_\mu}{\pi}.
  \label{eq:evidence}
\end{equation}
The workflow uses population annealing, because it can start from
$\beta=0$ and yields $\log(Z_\beta/Z_0)$ for the whole schedule in a
single run, whereas the ratios of replica exchange are referenced to
the lowest inverse temperature of the ladder.
The statistical properties of the population-annealing estimate of
$Z_\beta/Z_0$ are discussed in Refs.~\cite{Machta2010,Weigel2021};
its accuracy depends on the population size and on equilibration, as
illustrated in Section~\ref{sec:app:linreg}.

For this workflow, annealing starts from a uniform prior population
at $\beta=0$, either by placing $\beta=0$ at the beginning of the
schedule or through \code{anneal\_from\_beta0 = true}.
The output file \code{fx.txt} then contains $\log(Z_\beta/V_\Omega)$
directly, and no correction for the prior volume is needed (the
\code{-V} option of \code{odatse\_plt\_model\_evidence} keeps its
default).
A finite initial inverse temperature instead requires an additional
constant for the absolute evidence, derived in~\ref{app:normalization}.

The evidence maximum $\beta^\ast$ defines an empirical-Bayes estimate
of the effective residual noise, $\sigma_\mu^\ast=(2w_\mu\beta^\ast)^{-1/2}$.
The marginal posterior distributions of the parameters at
$\beta^\ast$ are conditional on that estimate; they do not marginalize
the uncertainty in $\beta$ or distinguish measurement noise from model
discrepancy.
If the evidence is used to compare different models, their
likelihoods must be normalized consistently and their priors stated
explicitly.

\begin{figure}[tb]
\centering
\includegraphics[width=\textwidth]{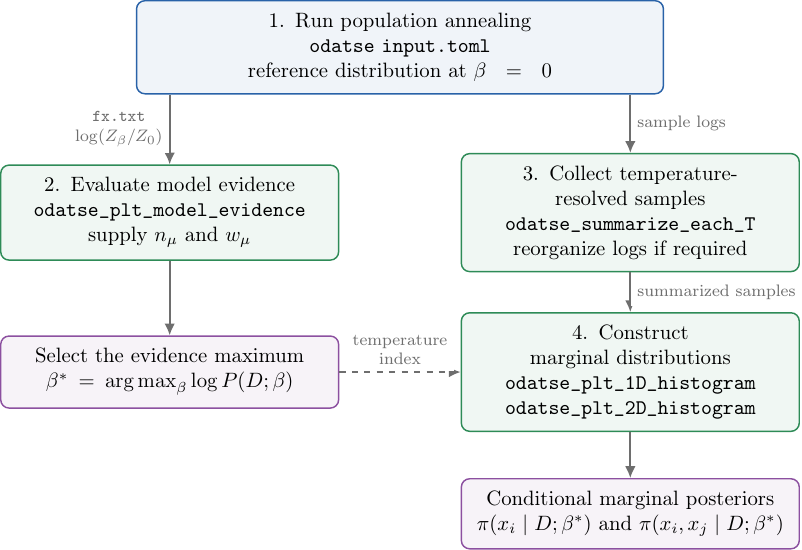}
\caption{Bayesian post-processing workflow.
Population annealing referenced to $\beta=0$ provides
partition-function ratios and sample logs.
The evidence calculation selects an effective noise
precision, while temperature-resolved samples provide
marginal posterior distributions.
The dashed arrow selects the temperature used for the
conditional posterior at $\beta^\ast$.
Blue denotes sampling, green denotes post-processing,
and purple denotes inferred quantities.}
\label{fig:workflow}
\end{figure}

The installed workflow (Figure~\ref{fig:workflow}) comprises:
\begin{enumerate}
\item Run \code{pamc} starting from $\beta=0$ as described above,
which yields the sample logs and \code{fx.txt}.
\item Use \code{odatse\_plt\_model\_evidence} with $n_\mu$ and
$w_\mu$ to evaluate Eq.~\eqref{eq:evidence} and locate its maximum
on the temperature schedule.
Independent runs provide a check of variability.
\item Use \code{odatse\_summarize\_each\_T} to collect the
samples at each temperature.
When necessary, expand combined output files beforehand with
\code{odatse\_extract\_combined}, or sort replica-exchange logs
by temperature with \code{odatse\_separateT}.
\item Apply \code{odatse\_plt\_1D\_histogram} and
\code{odatse\_plt\_2D\_histogram} to the samples at the selected
temperature to obtain the marginal posterior distributions of one or
two parameters.
\end{enumerate}

\subsection{Ecosystem: solver packages, documentation, and Gallery}
\label{sec:ecosystem}
\odatse{} is intended as a general platform for inverse-problem
analysis rather than as a tool for a particular experiment.
The application-specific direct problem solvers are therefore
separated from the core package and maintained independently, and
the core is complemented by documentation, templates, and worked
examples that allow users to connect their own simulation codes.
This subsection describes these components; their locations are given
in the Data availability statement.

\paragraph{Application-specific solver packages}
The three diffraction solvers bundled with \twodmat{}, for TRHEPD,
SXRD, and LEED, became independent packages in version~3.0 and are
maintained in the \code{2DMAT} organization on GitHub:
\code{odatse-STR} (TRHEPD) calls
\code{sim-trhepd-rheed}~\cite{Hanada2022},
\code{odatse-SXRD} (SXRD) calls \code{sxrdcalc}~\cite{Voegeli2006}, and
\code{odatse-LEED} (LEED) calls SATLEED~\cite{VanHove1993}.
A package for X-ray absorption fine structure, \code{odatse-XAFS},
which calls FEFF~\cite{Ankudinov1998}, has been added in the
\odatse{} generation.
All four are solvers of the external-program type described in
Section~\ref{sec:framework}: they prepare the input of the simulator,
run it, and compute the objective function from its output.
The external simulators are installed and licensed separately.
Beyond the structure analysis of surfaces, the same interface has been
applied to other fields of materials science:
\code{odatse-aenet} incorporates the machine-learning potentials of
{\ae}net~\cite{Artrith2016} as a surrogate for first-principles
calculations, and \code{odatse-specx} incorporates the
first-principles code AkaiKKR~\cite{Akai1989}.

\paragraph{Documentation, templates, and Gallery}
The user manual, including a chapter on customization, and the API
reference are part of the core package and are available online.
Tutorials and sample calculations are collected in the \odatse{}
Gallery, which contains the examples of Section~\ref{sec:applications}
such as the noise-level estimation by linear regression, the
Gaussian-process surrogate, and the QUBO problems, as well as sample
analyses for each solver package.
The Gallery also provides templates for writing solvers, from a
user-defined Python function optimized directly, through function
and solver modules that define a solver class, to a module that
calls an external program.
Through these libraries, documents, and examples, and through the
feedback of users who apply the framework to new problems, the
ecosystem around \odatse{} is being built up.

\section{Applications}
\label{sec:applications}

This section demonstrates the features of Section~\ref{sec:features}
on examples of increasing complexity.
Section~\ref{sec:app:linreg} validates the Bayesian workflow on a
one-parameter regression for which the evidence is known analytically;
Section~\ref{sec:app:ttopt} assesses tensor-train optimization on
benchmark functions; Section~\ref{sec:app:parallel} verifies the
two-level MPI parallelization with a solver whose work per evaluation
can be distributed; and Section~\ref{sec:app:trhepd} applies the
complete workflow to a structure analysis with synthetic TRHEPD data.
Further examples are given in~\ref{app:examples}: a fixed
surrogate model used as a solver, and binary and mixed-integer problems
solved with \code{ttopt}.

All calculations use \odatse{} v4.0.1 and were run on a Linux
workstation with two 32-core Intel Xeon Gold 6530 processors
(Python 3.13, NumPy 2.5, SciPy 1.18, Open MPI 5.0).
The versions of the solver package and of the external program used
in Section~\ref{sec:app:trhepd} are recorded in the data archive.
The inputs, selected outputs, and plotting scripts of all calculations
are archived in the ISSP Data Repository~\cite{Yoshimi2026data}, and
the examples are also provided as tutorials in the Gallery
(Section~\ref{sec:ecosystem}).

\subsection{Noise-level estimation by model evidence}
\label{sec:app:linreg}

A one-parameter regression $y=ax$ is an example for which the
partition function and the evidence are known in closed form; it
therefore serves to validate both the normalization estimate of
\code{pamc} and the evidence-based noise estimate.
The data are the six points $(x_i,y_i)$, $x_i=1,\ldots,6$, of the
Gallery example, generated from $y=x$ with additive noise.
A user-defined solver of about twenty lines, subclassing
\code{SolverBase} as described in Section~\ref{sec:framework}, returns
$f(a)=\sum_i(y_i-ax_i)^2$, which is Eq.~\eqref{eq:weighted_objective}
with $K=1$, $w_1=1$, and $n_1=6$; the prior is uniform on
$[-20,20]$.
Population annealing uses 100 logarithmically spaced inverse
temperatures from $10^{-5}$ to $10^2$ preceded by a $\beta=0$
reference step (\code{anneal\_from\_beta0 = true}),
100 Monte Carlo steps per temperature with proposal width 0.2,
and 1000 replicas on one process; execution takes about two minutes.

\begin{figure}[tb]
\centering
\includegraphics[width=\textwidth]{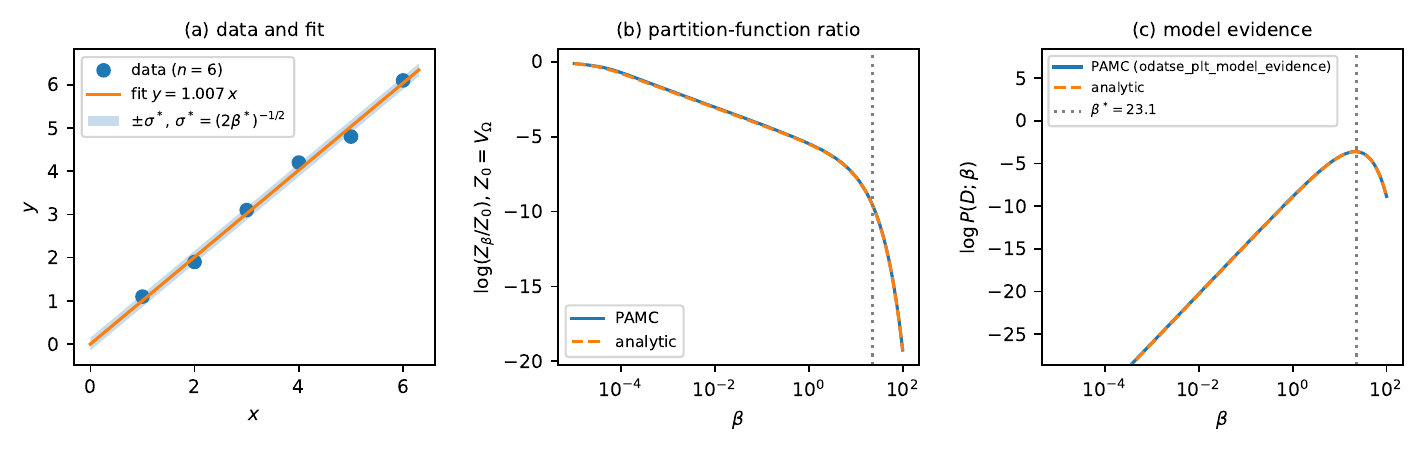}
\caption{Evidence-based noise estimation for linear regression.
(a) Synthetic data, the least-squares fit, and a band
of width $\pm\sigma^\ast$ around the fitted line.
This band indicates the estimated residual noise scale,
not a credible interval for the fitted slope.
(b) The partition-function ratio $\log(Z_\beta/Z_0)$
with $Z_0=V_\Omega$ from population annealing (solid)
and from Eq.~\eqref{eq:linreg_partition} (dashed).
(c) The log evidence computed with the functions of
\code{odatse\_plt\_model\_evidence} (solid) and its
analytical value (dashed).
The dotted lines mark the sampled maximum
$\beta^\ast=23.1$.}
\label{fig:linreg}
\end{figure}

Since the annealing starts from the uniform prior, \code{fx.txt}
reports $\log(Z_\beta/V_\Omega)$ and Eq.~\eqref{eq:evidence} gives
the absolute log evidence (Section~\ref{sec:bayes}).
Figure~\ref{fig:linreg}(b,c) compares the sampled curves with the
analytical expressions of~\ref{app:normalization}:
the partition-function ratio agrees to within 0.06 over the whole
schedule, and the two evidence curves coincide on the scale of the
figure.
We have also carried out a calculation with the default settings of
the Gallery example, which are chosen for a quick tutorial run:
10 replicas, a proposal width of 0.05, and a schedule starting at
$\beta=10^{-5}$ without the $\beta=0$ reference step; it takes a few
seconds.
Its partition-function ratio lies above the analytical value by
0.4--0.7 for $\beta>10^{-3}$, far more than the correction of 0.12
for the finite reference temperature derived in~\ref{app:normalization}.
The offset is caused by incomplete equilibration of the small
population during the first annealing steps; the analytical reference
makes such deficiencies visible.

The evidence of the 1000-replica run is maximal at the grid point
$\beta^\ast=23.1$, giving $\sigma^\ast=(2\beta^\ast)^{-1/2}=0.147$.
The analytical evidence has its maximum at $\beta=21.5$, of which
$\beta^\ast$ is the nearest point of the temperature grid (grid
spacing factor 1.18).
Because the posterior lies well within the prior bounds, this
maximum coincides with the Gaussian approximation
$\beta\simeq(n-1)/(2f^\ast)$, where $n=6$ is the number of data
points and $f^\ast$ the residual sum of squares of the least-squares
fit, so that the corresponding noise estimate
$\sqrt{f^\ast/(n-1)}=0.152$ is the residual standard error of the
fit; the difference from $\sigma^\ast$ reflects only the resolution
of the temperature grid.
This lightweight example thus tests the shape of the evidence curve
and the noise-scale estimate; as the comparison of the two runs
shows, the accuracy of the normalization estimate depends on the
population size and the proposal width, and in general also on the
number of Monte Carlo steps per temperature.

\subsection{Tensor-train optimization}
\label{sec:app:ttopt}

This example demonstrates the \code{ttopt} module and assesses its
performance on test functions provided by the \code{analytical}
solver of \odatse{}.
The suite consists of eleven functions that are commonly used to test
optimization methods~\cite{JamilYang2013,SurjanovicBingham2013} and,
except for the two-dimensional Himmelblau
and Rosenbrock functions, are defined for an arbitrary number of
variables; they are evaluated here in ten dimensions.
Most of them are multimodal with many local minima, so that the
suite tests the search over a large discretized space rather than
local refinement.

For each function we record the smallest objective value found
within a fixed budget of requested evaluations, together with the
wall-clock time and the fraction of requests served from the cache,
and compare the value with the known global minimum.
Table~\ref{tab:ttopt_bench} reports the results for $P_i=2$ and
$q_i=25$, that is, $2^{25}$ grid points per axis, $r_{\max}=4$, and
$10^5$ requested evaluations on eight MPI processes.
The suite takes approximately half a minute; the cache serves
52--82\% of the requests for the ten-dimensional functions and 95\%
for the two-dimensional Himmelblau function.
The grid does not contain the global minimizers of the functions,
because $2^{25}$ points on a symmetric interval exclude the origin,
so that the smallest attainable value is that of the grid point
nearest to the minimizer.
For Ackley, Exponential, Himmelblau, Qing, Rastrigin, and Schwefel, the search
reaches this grid point, and the difference from the reference
minimum is entirely due to the discretization (for Schwefel, also to
the rounded constant 418.9829 of the usual definition, which gives
$1.3\times10^{-4}$ at the minimizer); Alpine ends at another zero of
$x_i\sin x_i+0.1x_i$, with four coordinates at $x_i=-0.100$
($\sin x_i=-0.1$) and the others at the grid point nearest the
origin.
Griewank, Michalewicz, Rosenbrock, and Schaffer are not converged at
this budget; for Griewank, the search ends in a local minimum with
two coordinates at nonzero multiples of $\pi\sqrt{i}$.
The results illustrate representation and search capabilities
rather than a budget-matched ranking against other optimizers.

\begin{table}[tb]
\centering
\footnotesize
\caption{Analytical-function benchmarks with \code{ttopt}.
All calculations use $P_i=2$, $q_i=25$, $r_{\max}=4$,
a budget of $10^5$ requested evaluations, and eight
MPI processes.
The dimension is $d$.
Reference minima are those of the standard definitions, shown
to the displayed precision; the final column gives the smallest
objective value reported by the search, not a certified optimum.}
\label{tab:ttopt_bench}
\begin{tabular}{llrrr}
\toprule
Function & Search box & $d$ & Reference minimum & Best found \\
\midrule
Ackley & $[-32.768,32.768]^d$ & 10 & 0 & $3.9\times10^{-6}$ \\
Alpine & $[-10,10]^d$ & 10 & 0 & $2.9\times10^{-7}$ \\
Exponential & $[-1,1]^d$ & 10 & $-1$ & $-1.000$ \\
Griewank & $[-600,600]^d$ & 10 & 0 & $2.5\times10^{-2}$ \\
Himmelblau & $[-6,6]^d$ & 2 & 0 & $6.3\times10^{-13}$ \\
Michalewicz & $[0,\pi]^d$ & 10 & $-9.660$ & $-9.578$ \\
Qing & $[-500,500]^d$ & 10 & 0 & $1.5\times10^{-8}$ \\
Rastrigin & $[-5.12,5.12]^d$ & 10 & 0 & $4.6\times10^{-11}$ \\
Rosenbrock & $[-5,5]^d$ & 2 & 0 & $1.2\times10^{-2}$ \\
Schaffer & $[-100,100]^d$ & 10 & 0 & $3.6\times10^{-2}$ \\
Schwefel & $[-500,500]^d$ & 10 & 0 & $1.3\times10^{-4}$ \\
\bottomrule
\end{tabular}
\end{table}

Binary and mixed-variable problems handled by the same interface,
QUBO tests with up to 200 binary variables and three MINLPLib
problems with solver-defined penalties, are reported
in~\ref{app:discrete}.

\subsection{Two-level parallel execution}
\label{sec:app:parallel}

This example verifies that the two-level MPI parallelization of
Section~\ref{sec:mpi} works correctly, using a solver whose work per
evaluation is large enough for both levels to matter.
It is not a performance study: each configuration was run once.
The \code{parallel\_solver} sample searches, with \code{mapper}, for
the seed $s\in\{1,\ldots,20\}$ that minimizes the average of the
largest singular values of 50 random $1000\times1000$ matrices
generated from $s$.
\code{mapper} assigns contiguous blocks of seeds to the algorithm
ranks.
For each seed, the controller of a solver group generates the 50
matrices and broadcasts them over \code{solcomm} (about 400\,MB per
evaluation); every rank of the group computes the largest singular
values of its share of the matrices, the partial sums are combined by
an all-reduce, and the controller returns the average.
All runs use one BLAS thread per process (\code{OMP\_NUM\_THREADS=1})
so that each process occupies one core, with process binding disabled
(\code{-{}-bind-to none}).

\begin{figure}[tb]
\centering
\includegraphics[width=0.95\textwidth]{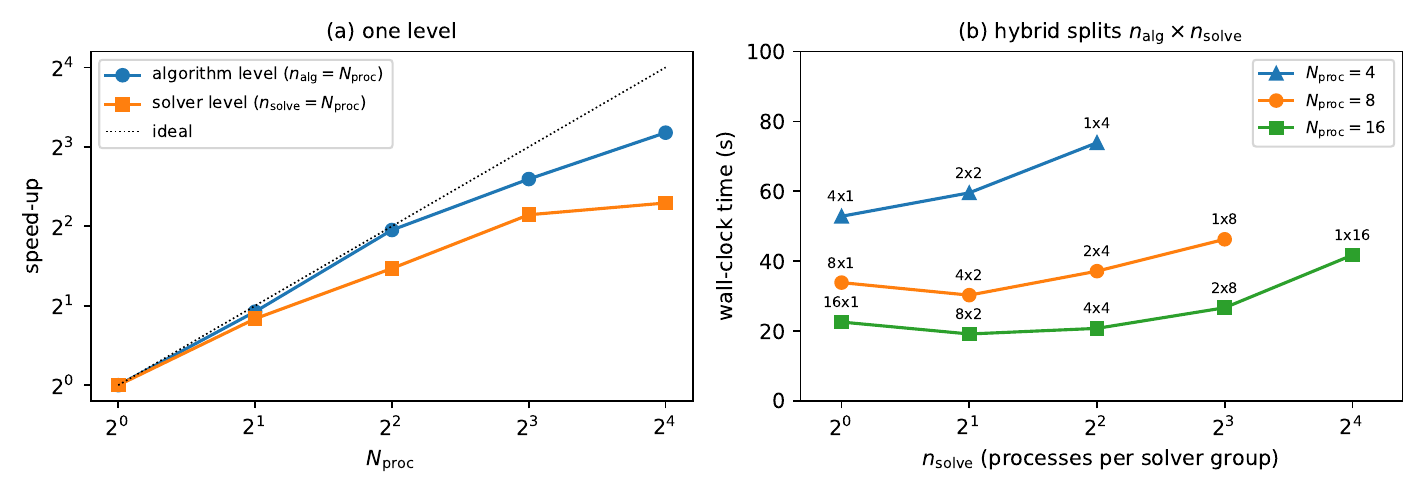}
\caption{Wall-clock times of the two-level parallel test on two
Intel Xeon Gold 6530 processors: 20 candidate seeds, each requiring
50 singular-value decompositions of $1000\times1000$ matrices, with
one thread per process and process binding disabled.
(a) Speed-up relative to the serial run with processes assigned only
to the algorithm level (circles) or only to the solver level
(squares); the dotted line indicates ideal linear scaling.
(b) Wall-clock times for decompositions
$n_{\mathrm{alg}}\times n_{\mathrm{solve}}$
using 4, 8, or 16 processes.
Each point is a single measurement of the whole \code{mpiexec} run,
not an average over repeated runs.}
\label{fig:scaling}
\end{figure}

All fifteen decompositions $n_{\mathrm{alg}}\times n_{\mathrm{solve}}$
with up to sixteen processes return the same minimum, $s=15$ with
$f=500.1006$, the values differing only in the last digits according
to the order of the floating-point reduction.
The wall-clock times (Figure~\ref{fig:scaling}) indicate that both
levels are effective: the serial run takes 204.7 seconds, and with
sixteen processes the mixed decompositions $8\times2$ and $4\times4$
take 19.2 and 20.8 seconds, compared with 22.6 seconds for
$16\times1$ and 41.8 seconds for $1\times16$.
Solver-level parallelism alone scales worse because the matrices are
generated by the controller and broadcast, and because 50
decompositions do not divide evenly over sixteen ranks; at the
algorithm level, 20 seeds do not divide evenly over eight or sixteen
groups.
The best split therefore depends on the workload, which is why it is
left to the user.

\subsection{Bayesian structure analysis using simulated TRHEPD data}
\label{sec:app:trhepd}

The last example takes up the TRHEPD analysis of the
Ge(001)-c($4\times2$) surface from the \twodmat{}
paper~\cite{Motoyama2022} and carries it through the Bayesian
workflow of Section~\ref{sec:bayes}.
There, artificial noise-free data generated from a reference
structure were analyzed by local optimization, grid search, Bayesian
optimization, and replica exchange, which revealed two minima related
by the exchange of two coordinates and a valley of the objective
function along their difference.
Here we add noise to the data, estimate the noise level from the
evidence, and obtain the posterior distribution at the selected
noise precision, which was not possible with \twodmat{} because it did
not provide partition functions.
Lengths in this example are given in the non-SI unit \angstrom{}
($1\,\angstrom=0.1\,\text{nm}$), which is customary in surface
crystallography.

The parameters $z_1$, $z_2$, and $z_3$ are the heights of the first
three atomic layers of the reconstructed surface (Fig.~3 of
Ref.~\cite{Motoyama2022}); the deeper layers are fixed.
The data are a rocking curve, the diffraction intensity at 70
glancing angles from $0.1^\circ$ to $7.0^\circ$, under the one-beam
condition~\cite{Ichimiya1987}, in which the intensity is sensitive
mainly to the layer heights.
The solver is \code{odatse-STR}, which runs the external simulator
\code{sim-trhepd-rheed}~\cite{Hanada2022} for every evaluation
(Section~\ref{sec:ecosystem}).
The reference curve is generated by the same simulator from
$(z_1,z_2,z_3)=(5.231,4.371,3.596)\,\angstrom$~\cite{Shirasawa2006},
and the prior is uniform over $[3,6]^3\,\angstrom^3$.
Each of the 70 intensities is multiplied by $1+\epsilon_i$, with
independent Gaussian $\epsilon_i$ of standard deviation 0.1, and the
curve is normalized to unit total intensity; the realized
root-mean-square perturbation after normalization is
$1.36\times10^{-3}$.
The solver normalizes its prediction in the same way and returns
\begin{equation}
 f(z_1,z_2,z_3)=
 \sum_{i=1}^{70}(I_{\mathrm{cal},i}-I_{\mathrm{ref},i})^2,
\end{equation}
selected by \code{Rfactor\_type = "A2"}; the default \code{"A"}
returns its square root, the R-factor used in Ref.~\cite{Motoyama2022}.
Inference uses the independent additive Gaussian likelihood of
Eq.~\eqref{eq:likelihood}; it is a working approximation, because
the synthetic noise is multiplicative and the normalization
correlates the intensities.

Population annealing uses 800 replicas on eight processes,
51 logarithmically spaced inverse temperatures from 1 to $10^6$,
and \code{anneal\_from\_beta0 = true}.
Forty Monte Carlo steps per temperature with proposal width
$0.01\,\angstrom$ give approximately $1.6\times10^6$ solver calls,
that is, runs of the external simulator, in 26--28 minutes.
Two sampling seeds are used.
The proposal width retains an acceptance rate of about 0.45 at the
selected temperature; in a preliminary run with width
$0.1\,\angstrom$, the acceptance rate dropped to nearly zero at low
temperatures and fewer than a quarter of the walkers remained
distinct at the lowest temperature.

\begin{figure}[tb]
\centering
\includegraphics[width=\textwidth]{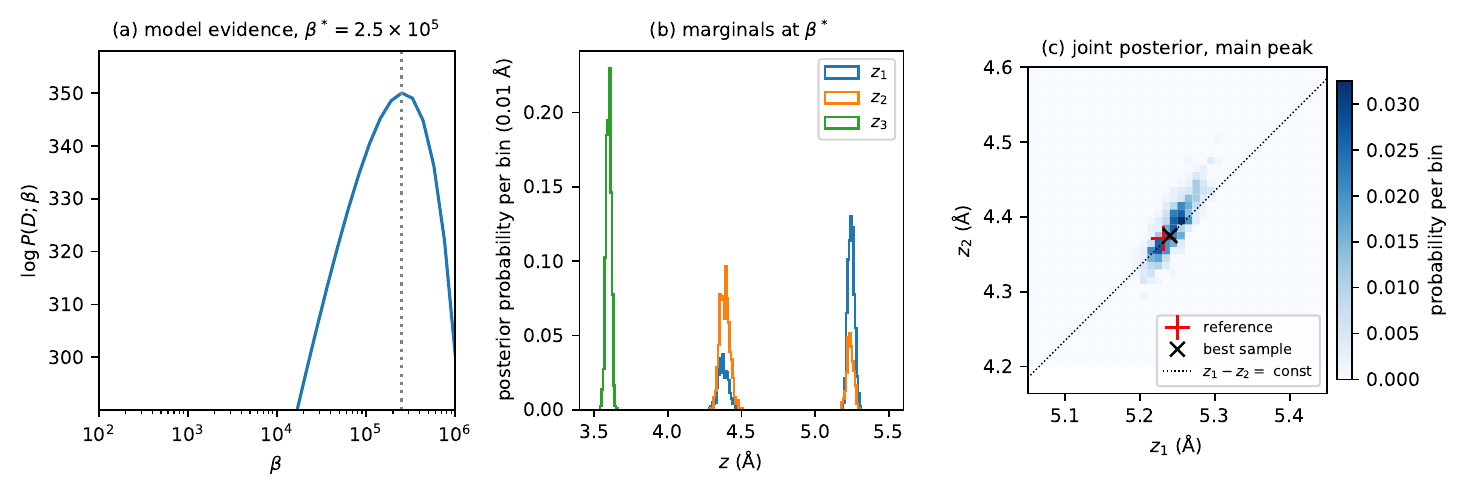}
\caption{Bayesian analysis of synthetic TRHEPD data
for Ge(001)-c($4\times2$) using population annealing
and \code{odatse-STR}.
(a) Estimated log evidence under the working
independent Gaussian likelihood.
The dotted line identifies the temperature-grid
maximum $\beta^\ast$.
Evidence here is a density for continuous data,
so its logarithm need not be negative.
(b) Marginal distributions at $\beta^\ast$
for the second sampling seed, shown as probability
per $0.01\,\angstrom$ bin.
The two modes in $z_1$ and $z_2$ reflect
exchange symmetry; their unequal sampled weights
indicate incomplete global equilibration.
(c) Joint distribution restricted to the mode
with $z_1>z_2$, also for the second seed.
The red plus marks the reference structure,
the black cross marks the best sampled point,
and the dotted line indicates
$z_1-z_2=\mathrm{const}$.
The plotted distributions are conditional on
the selected noise precision.}
\label{fig:evidence}
\end{figure}

Both runs select $\beta^\ast=2.5\times10^5$ on the temperature grid
(spacing factor 1.32; Figure~\ref{fig:evidence}(a)).
The local Laplace estimate $(n-N)/(2f^\ast)$, with $n=70$ data
points, $N=3$ parameters, and the smallest sampled value
$f^\ast=1.28\times10^{-4}$, gives $2.6\times10^5$, consistent with
this value.
The effective noise scale $\sigma^\ast=(2\beta^\ast)^{-1/2}
=1.41\times10^{-3}$ is close to the imposed perturbation scale.
This agreement concerns the aggregate residual noise; it does not
establish calibrated intervals under the correlated noise process.

The marginal distributions at $\beta^\ast$ (Figure~\ref{fig:evidence}(b))
show a single peak of $z_3$ near $3.60\,\angstrom$ with standard
deviation $0.017\,\angstrom$ and, for $z_1$ and $z_2$, the two modes
near $(5.24,4.38)$ and $(4.38,5.24)\,\angstrom$ related by the
exchange symmetry.
Their true weights are equal, but the sampled populations split
approximately 81:19 and 73:27 in the two runs: both solutions are
recovered, but the global mode weights are not converged, in
contrast to the symmetric histograms obtained for noise-free data in
Ref.~\cite{Motoyama2022}.
Within the mode with $z_1>z_2$ (Figure~\ref{fig:evidence}(c)), the
distribution is elongated along $z_1-z_2=\mathrm{const}$, as the
objective function is most sensitive to the distance between the
first two layers~\cite{Motoyama2022,Hoshi2022}: the standard
deviations along $(z_1+z_2)/\sqrt{2}$ and $(z_1-z_2)/\sqrt{2}$ are
0.037 and $0.013\,\angstrom$, and those of $z_1$ and $z_2$ are 0.021
and $0.033\,\angstrom$, with agreement between the two runs within
$0.003\,\angstrom$.
The best samples, $(5.242,4.379,3.596)$ and
$(5.240,4.375,3.594)\,\angstrom$, are close to the reference.

These results connect the evidence-based noise estimate with the
structural correlations that a single optimum cannot describe.
The widths characterize a local mode conditional on $\beta^\ast$ and
on the chosen likelihood; the convergence of the global weights and
of the absolute evidence requires further sampling tests.

\section{Summary}
\label{sec:summary}

We have described \odatse{}, the successor of \twodmat{}, in its
current version~4.0.
The reorganization separated the search and sampling
algorithms from the application-specific solvers, which are now
maintained as independent packages, and changed the license to
MPL-2.0, which eases the use of the framework as a library in other
software.
On this basis, versions~3 and~4 add a common checkpoint and
restart interface, random search, two-level MPI parallelization, the global optimizers
of SciPy, tensor-train optimization, a choice of local optimization
methods, and a Bayesian workflow from partition-function estimates to
evidence-based noise estimation and posterior distributions.

The examples exercise these features one by one.
The linear-regression example validates the partition-function
estimate and the noise estimate against closed-form results and shows
that the analytical reference exposes insufficient equilibration.
The tensor-train benchmark reaches the grid resolution for most test
functions within a fixed evaluation budget.
The two-level parallel test returns the same optimum for all process
decompositions.
The TRHEPD example carries the analysis of the \twodmat{} paper to a
Bayesian analysis: it recovers the two symmetric structural solutions
and their correlation, and estimates a noise scale close to the
imposed perturbation, while the unequal weights of the two modes show
that further sampling is needed for global convergence.

Some limitations remain.
The SciPy optimizers of \code{minsearch} and \code{global\_search}
cannot be restarted.
Future work includes adaptive surrogate models and the extension to
further classes of problems by connecting other existing simulation
codes as solvers; through the Gallery and the feedback of users who
apply the framework to new problems, \odatse{} is being developed as
a general platform for inverse-problem analysis.

\section*{Acknowledgements}
\twodmat{} and \odatse{} have been developed and maintained under the ``Project for advancement of software usability in materials science'' (PASUMS)~\cite{Yoshimi2025} of the Institute for Solid State Physics, The University of Tokyo, which supported the development in FY2020, FY2021, and FY2024.
The development was also supported by JSPS KAKENHI Grant Numbers 19H04125, 20H00581, 21K19773, 23K18465, and 25K15153.
All authors were supported by the Japan Science and Technology Agency (JST) through the Moonshot Research and Development Program, Grant Number JPMJMS24A3, for the development of version 4, the applications reported in Section~\ref{sec:applications}, and the preparation of this manuscript.
The program code was developed on the supercomputer Fugaku through the HPCI System Research Project (Project ID: hp220248, hp230304, hp240304, hp250295).
We thank K. Akamatsu, K. Tsukamoto, and H. Shimizu for their contributions to the code, the manual, and the tutorials.

\section*{Declaration of competing interests}
The authors declare that they have no known competing financial
interests or personal relationships that could have appeared to
influence the work reported in this paper.

\section*{Data availability}

\odatse{} v4.0.1 is available under MPL-2.0 from
\url{https://github.com/issp-center-dev/ODAT-SE} and PyPI.
The manual and API reference are at
\url{https://issp-center-dev.github.io/ODAT-SE/manual/main/en/}.
The diffraction and spectroscopy solver packages
(\code{odatse-STR}, \code{odatse-SXRD}, \code{odatse-LEED},
and \code{odatse-XAFS}) are hosted under
\url{https://github.com/2DMAT}.
Related interfaces are available at
\url{https://github.com/k-yoshimi/odatse-aenet} and
\url{https://github.com/k-yoshimi/odatse-specx}.
Tutorials, sample inputs, and the solver templates are in the Gallery:
\url{https://isspns-gitlab.issp.u-tokyo.ac.jp/takeohoshi/odat-se-gallery}
(templates under \code{data/tutorial/solver-template}).
Inputs, selected outputs, and plotting scripts for the calculations
are archived in the ISSP Data Repository~\cite{Yoshimi2026data}.
Large per-evaluation simulator logs and full TTOpt histories are
excluded; the archived inputs specify their regeneration.

\appendix

\section{Evidence normalization}
\label{app:normalization}

\subsection{Finite reference temperature}
The samplers of Section~\ref{sec:mcext} estimate partition-function
ratios relative to the first inverse temperature of the schedule,
$\beta_{\mathrm{ref}}$: this is the case for population annealing
without \code{anneal\_from\_beta0} and for replica exchange, whose
ratios are referenced to the lowest inverse temperature of the ladder.
The absolute evidence of Eq.~\eqref{eq:evidence} then requires an
additional constant,
\begin{equation}
 \begin{split}
 \log P(D;\beta)={}&
 \log\frac{Z_\beta}{Z_{\beta_{\mathrm{ref}}}}
 +\log\frac{Z_{\beta_{\mathrm{ref}}}}{V_\Omega}\\
 &+\sum_\mu\frac{n_\mu}{2}\log\frac{\beta w_\mu}{\pi}.
 \end{split}
 \label{eq:evidence_finite_reference}
\end{equation}
The second term is independent of $\beta$, so that the location of
the evidence maximum $\beta^\ast$ does not depend on it; it matters
only for the absolute value of the evidence and hence for the
comparison of different models.
It vanishes for $\beta_{\mathrm{ref}}=0$ with a uniform prior, which
is the case for the runs of Sections~\ref{sec:app:linreg}
and~\ref{sec:app:trhepd}; otherwise $Z_{\beta_{\mathrm{ref}}}$ has to
be obtained separately, for example from a short run from $\beta=0$
to $\beta_{\mathrm{ref}}$ or, as below, analytically.
For Monte Carlo sampling on a discrete mesh (Section~\ref{sec:mcext})
with a uniform prior, the integrals become sums and $V_\Omega$ is the
number of admissible mesh points.

\subsection{Partition function of the linear regression}
For the regression $y=ax$ of Section~\ref{sec:app:linreg} with data
$(x_i,y_i)$, $i=1,\ldots,n$, let $S=\sum_i x_i^2$,
$\widehat a=\sum_i x_i y_i/S$, and $f^\ast=f(\widehat a)$.
Then $f(a)=f^\ast+S(a-\widehat a)^2$, and for a uniform prior on
$[a_-,a_+]$ and $\beta>0$,
\begin{equation}
 \begin{split}
 Z_\beta={}&e^{-\beta f^\ast}\frac{\sqrt{\pi}}{2\sqrt{\beta S}}
 \biggl[
 \operatorname{erf}\!\left(\sqrt{\beta S}(a_+-\widehat a)\right)\\
 &\hspace{1.2cm}
 -\operatorname{erf}\!\left(\sqrt{\beta S}(a_--\widehat a)\right)
 \biggr],
 \end{split}
 \label{eq:linreg_partition}
\end{equation}
while $Z_0=V_\Omega=a_+-a_-$.
The evidence follows from Eq.~\eqref{eq:evidence} with $K=1$,
$w_1=1$, and $n_1=n$.
For the archived data, $n=6$, $x_i=i$, $S=91$, $\widehat a=1.0066$,
$f^\ast=0.1160$, $a_\mp=\mp20$, and $V_\Omega=40$; these values give
the dashed curves in Figure~\ref{fig:linreg}(b,c).
For the run with the Gallery default settings in
Section~\ref{sec:app:linreg}, which starts at
$\beta_{\mathrm{ref}}=10^{-5}$, the correction
$\log(V_\Omega/Z_{\beta_{\mathrm{ref}}})$ evaluates to only 0.12, so
that the offset of 0.4--0.7 observed there is due to the sampling,
not to the normalization.

\section{Additional solver and optimization examples}
\label{app:examples}

\subsection{Fixed-surrogate integration}
\label{app:surrogate}

This example shows that a machine-learning model can be used as a
solver through the same interface as a physical simulation, which is
the pattern of \code{odatse-aenet} (Section~\ref{sec:ecosystem}) in a
minimal form.
A Gaussian-process regressor~\cite{RasmussenWilliams2005} is trained
on 500 evaluations of the Himmelblau function at uniformly distributed
points in $[-6,6]^2$, with Gaussian noise of standard deviation 0.01
added to the training values; inputs and targets are standardized for
training.
The kernel is a squared-exponential kernel with a white-noise term,
whose hyperparameters are fitted by maximizing the marginal likelihood
with scikit-learn~\cite{Pedregosa2011} 1.9, an example-specific
dependency.
The predictor is wrapped as an \odatse{} solver by a
\code{SolverBase} subclass of a few lines whose \code{evaluate}
returns the predicted mean; the solver is constructed in a Python
script that drives the calculation (Section~\ref{sec:whatis}), and
the \code{[solver]} section of the input file is left empty.
The random seed of the input file is used both for the data and for
the fit.

\begin{figure}[tb]
\centering
\includegraphics[width=\textwidth]{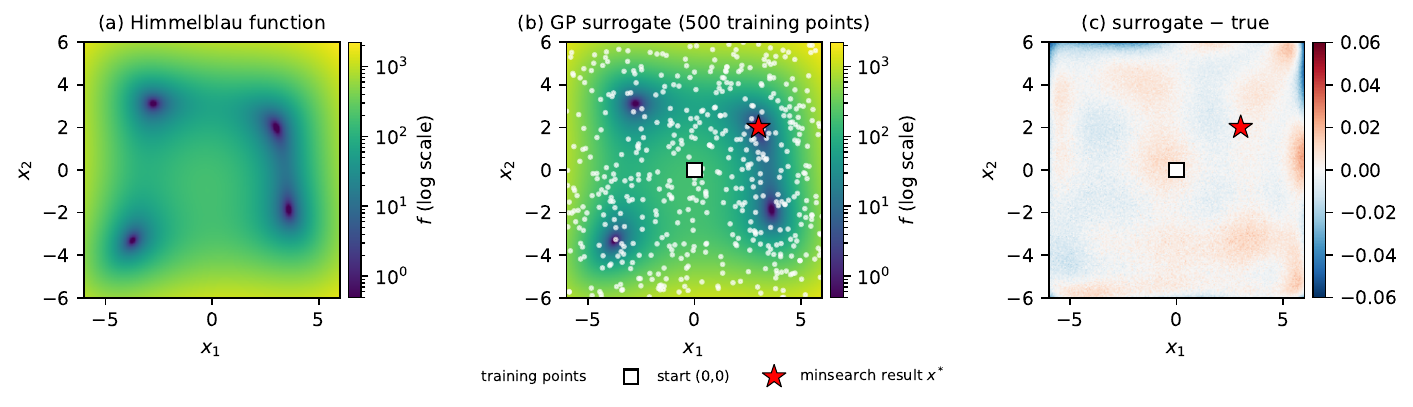}
\caption{Optimization using a fixed Gaussian-process surrogate.
(a) The Himmelblau function.
(b) The surrogate used as the \odatse{} solver, with
the 500 training locations shown as dots.
The square marks the initial point of \code{minsearch},
and the star marks the minimum found.
(c) The signed difference between the surrogate and
the original function.
Panels (a) and (b) use logarithmic color scales;
panel (c) uses a linear scale.
The surrogate is trained before optimization and is
not updated during the search.}
\label{fig:gpr}
\end{figure}

On 500 held-out points, the mean squared error of the surrogate
against the noise-free function is $4\times10^{-5}$, that is, a
root-mean-square error of about $7\times10^{-3}$ in the original
function units; Figure~\ref{fig:gpr}(c) shows that the signed error
stays within $\pm0.07$ over the box, with the largest deviations
at its boundary.
Starting from the origin, \code{minsearch} on the surrogate reaches
$(3.006,1.999)$, close to the minimum $(3,2)$ of the original
function; the surrogate value there is $-2\times10^{-3}$, slightly
below the true minimum 0, and the original function evaluates to
$1\times10^{-3}$ at the optimized point.
These numbers depend on the computing environment, because the
hyperparameter optimization of the Gaussian process is sensitive to
differences in floating-point arithmetic: a rerun with the same seed
on another machine gave a mean squared error of $5\times10^{-5}$ and
the optimum $(3.002,1.996)$, with the same qualitative result.
The training set is fixed: this is an interface example, not an
adaptive surrogate-assisted algorithm.
When such a predictor is used as the objective of a sampler, the
uncertainty of the Gaussian process is not propagated automatically,
and final candidates should be verified with the original forward
model.

\subsection{Binary and mixed-integer problems}
\label{app:discrete}

This appendix supplements the continuous benchmarks of
Section~\ref{sec:app:ttopt} with binary and mixed-variable problems.
In \code{ttopt}, a binary variable is a single index with $P_i=2$ and
$q_i=1$, and an integer variable with $m$ values a single index with
$P_i=m$ (Section~\ref{sec:newalg}).

Quadratic unconstrained binary optimization (QUBO) minimizes
$x^{\mathsf T}Qx$ over $x\in\{0,1\}^N$~\cite{Lucas2014}.
Random instances generated with \code{qubogen} 0.1.1 encode max-cut
and minimum vertex cover on random graphs with edge probability 0.5,
and quadratic knapsack with random values and weights, for
$N=50,100,200$.
Ten instances per class and size are solved with ten optimization
seeds each, using $r_{\max}=4$, $10^5$ requested evaluations per
seed, and sixteen MPI processes.
As a reference, the Gallery distributes the solutions of the same
instances obtained with the tabu-search QUBO solver of MATLAB with
default settings; these are heuristic solutions, not certified
optima.
For each instance, the best of the ten runs and each single run are
compared with the reference through the relative deviation
$|f-f_{\mathrm{ref}}|/|f_{\mathrm{ref}}|$ of the QUBO values.
Table~\ref{tab:qubo} gives the mean and maximum over the ten
instances of the best-of-ten deviation, the mean deviation of the
hundred single runs, the number of instances whose best-of-ten
value matches the reference within $10^{-5}$, and the mean wall-clock
time per run, 4--13 seconds.
All quadratic-knapsack instances with $N=50$ and $100$ and most with
$N=200$ match the reference; the minimum-vertex-cover deviations are at
most 0.03\%, and the max-cut deviations grow with $N$ up to 0.47\%.
The deviations refer to the encoded QUBO values, which contain large
offsets and penalty terms; for minimum vertex cover, for example, a
cover larger by one vertex changes the value by only about
$10^{-4}$ relative, which is why instances with a displayed deviation
of 0.00\% need not match the reference.

\begin{table}[tb]
\centering
\footnotesize
\caption{QUBO results with \code{ttopt}, using $10^5$
requested evaluations per seed and 16 MPI processes.
Best-of-ten columns give the mean and maximum
relative deviations from the tabu-search reference
over ten instances.
``Single'' is the mean deviation over the hundred
individual runs.
``Matches'' counts instances whose best-of-ten
objective agrees with the reference to a relative
tolerance of $10^{-5}$; it does not certify
global optimality.
Times are mean wall-clock times per run.
A displayed deviation of 0.00\% can be a rounded
nonzero value.}
\label{tab:qubo}
\begin{tabular}{lrrrrcr}
\toprule
Problem & $N$ & \multicolumn{2}{c}{Best of 10 (\%)} & Single (\%) & Matches & Time (s) \\
 & & Mean & Max & Mean & & \\
\midrule
Max-cut & 50 & 0.03 & 0.28 & 0.74 & 9/10 & 4 \\
Max-cut & 100 & 0.07 & 0.21 & 0.74 & 3/10 & 7 \\
Max-cut & 200 & 0.31 & 0.47 & 0.68 & 0/10 & 13 \\
Minimum vertex cover & 50 & 0.01 & 0.03 & 0.02 & 5/10 & 4 \\
Minimum vertex cover & 100 & 0.00 & 0.00 & 0.01 & 3/10 & 7 \\
Minimum vertex cover & 200 & 0.00 & 0.00 & 0.00 & 0/10 & 12 \\
Quadratic knapsack & 50 & 0.00 & 0.00 & 0.06 & 10/10 & 4 \\
Quadratic knapsack & 100 & 0.00 & 0.00 & 0.37 & 10/10 & 6 \\
Quadratic knapsack & 200 & 0.02 & 0.14 & 0.42 & 8/10 & 12 \\
\bottomrule
\end{tabular}
\end{table}

Three MINLPLib problems~\cite{Bussieck2003} test mixed variables:
\code{alan}, \code{synthes2}, and \code{cvxnonsep\_normcon20}
(Table~\ref{tab:mip}).
Equality constraints are eliminated by substitution, which leaves 2
continuous and 4 binary variables for \code{alan} and 6 continuous
and 4 binary variables for \code{synthes2}; \code{cvxnonsep\_normcon20}
has 10 continuous and 10 integer variables in $\{0,\ldots,5\}$ and a
single norm constraint.
Continuous axes use $P_i=2$ and $q_i=20$, binary axes $P_i=2$ and
$q_i=1$, and integer axes $P_i=6$ and $q_i=1$.
The solver returns $g(x)+\lambda\sum_k\max(0,c_k(x))^2$, where $g$ is
the objective and $c_k(x)\le0$ are the remaining inequality
constraints (for \code{synthes2}, four inequalities implied by the
others are omitted).
The penalty weight $\lambda$ starts at 10, 100, and 1 for the three
problems, respectively, and doubles over 20 stages; the best points of
the previous stages are passed to the next stage through
\code{init\_points}.
Each stage uses $10^5$ requested evaluations on eight processes.

Table~\ref{tab:mip} compares the objective at the final point,
without the penalty term, with the global optimum reported by
MINLPLib.
The result for \code{alan} agrees to the displayed precision, that for
\code{synthes2} differs by $4\times10^{-4}$\,\%, and that for
\code{cvxnonsep\_normcon20} by 0.46\% of the optimum, where the search
ends at an integer assignment different from the reference solution.
The constraint violation $\sum_k\max(0,c_k)^2$ at the final points is
of order $10^{-12}$ or zero, which indicates numerical feasibility,
not an exact certificate.

\begin{table}[tb]
\centering
\footnotesize
\caption{Mixed-integer programs treated with \code{ttopt}
and a staged quadratic penalty.
The numbers of variables (c: continuous, b: binary, i: integer) and
of constraints (eq: equality, ineq: inequality) are those of the
MINLPLib formulation; the equality constraints are eliminated by
substitution, so that the search space has fewer variables (see
text).
Each stage uses eight MPI processes and a budget of $10^5$ requested
evaluations.
``Best found'' is the objective at the final point without the
penalty term, ``Reference'' the global optimum reported by MINLPLib
to the displayed precision, and ``Violation'' is
$\sum_k\max(0,c_k)^2$ at the final point.}
\label{tab:mip}
\setlength{\tabcolsep}{4pt}%
\begin{tabular}{@{}llrrrr@{}}
\toprule
Instance & Variables & Constraints & Best found & Reference & Violation \\
\midrule
\code{alan} & 4 c, 4 b & 2 eq, 5 ineq & 2.9250 & 2.9250 & 0 \\
\code{synthes2} & 6 c, 5 b & 1 eq, 13 ineq & 73.036 & 73.035 & $5\times10^{-12}$ \\
\code{cvxnonsep\_normcon20} & 10 c, 10 i & 0 eq, 1 ineq & $-21.650$ & $-21.749$ & $6\times10^{-12}$ \\
\bottomrule
\end{tabular}
\end{table}

\bibliographystyle{elsarticle-num}
\bibliography{odatse}

\begin{thebibliography}{10}
\expandafter\ifx\csname url\endcsname\relax
  \def\url#1{\texttt{#1}}\fi
\expandafter\ifx\csname urlprefix\endcsname\relax\def\urlprefix{URL }\fi
\expandafter\ifx\csname href\endcsname\relax
  \def\href#1#2{#2} \def\path#1{#1}\fi

\bibitem{Nagata2012}
K.~Nagata, S.~Sugita, M.~Okada, Bayesian spectral deconvolution with the
  exchange {Monte Carlo} method, Neural Networks 28 (2012) 82--89.
\newblock \href {https://doi.org/10.1016/j.neunet.2011.12.001}
  {\path{doi:10.1016/j.neunet.2011.12.001}}.

\bibitem{Hugenschmidt2016}
C.~Hugenschmidt, Positrons in surface physics, Surface Science Reports 71
  (2016) 547--594.
\newblock \href {https://doi.org/10.1016/j.surfrep.2016.09.002}
  {\path{doi:10.1016/j.surfrep.2016.09.002}}.

\bibitem{Fukaya2019}
Y.~Fukaya, A.~Kawasuso, A.~Ichimiya, T.~Hyodo, Total-reflection high-energy
  positron diffraction ({TRHEPD}) for structure determination of the topmost
  and immediate sub-surface atomic layers, Journal of Physics D: Applied
  Physics 52 (2019) 013002.
\newblock \href {https://doi.org/10.1088/1361-6463/aadf14}
  {\path{doi:10.1088/1361-6463/aadf14}}.

\bibitem{Feidenhansl1989}
R.~Feidenhans'l, Surface structure determination by x-ray diffraction, Surface
  Science Reports 10 (1989) 105--188.
\newblock \href {https://doi.org/10.1016/0167-5729(89)90002-2}
  {\path{doi:10.1016/0167-5729(89)90002-2}}.

\bibitem{Tajiri2020}
H.~Tajiri, Progress in surface x-ray crystallography and the phase problem,
  Japanese Journal of Applied Physics 59 (2020) 020503.
\newblock \href {https://doi.org/10.7567/1347-4065/ab631e}
  {\path{doi:10.7567/1347-4065/ab631e}}.

\bibitem{VanHove1993}
M.~A. Van~Hove, W.~Moritz, H.~Over, P.~J. Rous, A.~Wander, A.~Barbieri,
  N.~Materer, U.~Starke, G.~A. Somorjai, Automated determination of complex
  surface structures by {LEED}, Surface Science Reports 19 (1993) 191--229.
\newblock \href {https://doi.org/10.1016/0167-5729(93)90011-D}
  {\path{doi:10.1016/0167-5729(93)90011-D}}.

\bibitem{Hoshi2022}
T.~Hoshi, D.~Sakata, S.~Oie, I.~Mochizuki, S.~Tanaka, T.~Hyodo, K.~Hukushima,
  Data-driven sensitivity analysis in surface structure determination using
  total-reflection high-energy positron diffraction ({TRHEPD}), Computer
  Physics Communications 271 (2022) 108186.
\newblock \href {https://doi.org/10.1016/j.cpc.2021.108186}
  {\path{doi:10.1016/j.cpc.2021.108186}}.

\bibitem{Anada2017}
M.~Anada, Y.~Nakanishi-Ohno, M.~Okada, T.~Kimura, Y.~Wakabayashi, Bayesian
  inference of metal oxide ultrathin film structure based on crystal truncation
  rod measurements, Journal of Applied Crystallography 50 (2017) 1611--1616.
\newblock \href {https://doi.org/10.1107/S1600576717013292}
  {\path{doi:10.1107/S1600576717013292}}.

\bibitem{Nagai2020}
K.~Nagai, M.~Anada, Y.~Nakanishi-Ohno, M.~Okada, Y.~Wakabayashi, Robust surface
  structure analysis with reliable uncertainty estimation using the exchange
  {Monte Carlo} method, Journal of Applied Crystallography 53 (2020) 387--392.
\newblock \href {https://doi.org/10.1107/S1600576720001314}
  {\path{doi:10.1107/S1600576720001314}}.

\bibitem{Motoyama2022}
Y.~Motoyama, K.~Yoshimi, I.~Mochizuki, H.~Iwamoto, H.~Ichinose, T.~Hoshi,
  Data-analysis software framework {2DMAT} and its application to experimental
  measurements for two-dimensional material structures, Computer Physics
  Communications 280 (2022) 108465.
\newblock \href {https://doi.org/10.1016/j.cpc.2022.108465}
  {\path{doi:10.1016/j.cpc.2022.108465}}.

\bibitem{Tanaka2020}
K.~Tanaka, T.~Hoshi, I.~Mochizuki, T.~Hanada, A.~Ichimiya, T.~Hyodo,
  Development of data-analysis software for total-reflection high-energy
  positron diffraction ({TRHEPD}), Acta Physica Polonica A 137 (2020) 188--192.
\newblock \href {https://doi.org/10.12693/APhysPolA.137.188}
  {\path{doi:10.12693/APhysPolA.137.188}}.

\bibitem{Tanaka2023}
K.~Tanaka, I.~Mochizuki, T.~Hanada, A.~Ichimiya, T.~Hyodo, T.~Hoshi, A
  two-stage data-analysis method for total-reflection high-energy positron
  diffraction ({TRHEPD}), JJAP Conference Proceedings 9 (2023) 011301.
\newblock \href {https://doi.org/10.56646/jjapcp.9.0_011301}
  {\path{doi:10.56646/jjapcp.9.0_011301}}.

\bibitem{Tsujikawa2022}
Y.~Tsujikawa, M.~Shoji, M.~Hamada, T.~Takeda, I.~Mochizuki, T.~Hyodo,
  I.~Matsuda, A.~Takayama, Structure of $\chi_3$-borophene studied by
  total-reflection high-energy positron diffraction ({TRHEPD}), Molecules 27
  (2022) 4219.
\newblock \href {https://doi.org/10.3390/molecules27134219}
  {\path{doi:10.3390/molecules27134219}}.

\bibitem{Sato2025}
Y.~Sato, Y.~Fukaya, A.~Nakano, T.~Hoshi, C.-C. Lee, K.~Yoshimi, T.~Ozaki,
  T.~Nakashima, Y.~Ando, H.~Aoyama, T.~Abukawa, Y.~Tsujikawa, M.~Horio,
  M.~Niibe, F.~Komori, I.~Matsuda, Surface structure of the $3\times3$-{Si}
  phase on {Al(111)}, studied by the multiple usages of positron diffraction
  and core-level photoemission spectroscopy, Physical Review Materials 9 (2025)
  014002.
\newblock \href {https://doi.org/10.1103/PhysRevMaterials.9.014002}
  {\path{doi:10.1103/PhysRevMaterials.9.014002}}.

\bibitem{Liu2025}
S.~Liu, M.~Roppongi, M.~Kimata, K.~Ishihara, R.~Grasset, M.~Konczykowski, B.~R.
  Ortiz, S.~D. Wilson, K.~Yoshimi, T.~Shibauchi, K.~Hashimoto, Impact of tiny
  {Fermi} pockets with extremely high mobility on the {Hall} anomaly in the
  kagome metal {CsV$_3$Sb$_5$}, Physical Review Letters 135 (2025) 056502.
\newblock \href {https://doi.org/10.1103/d4dw-2v6k}
  {\path{doi:10.1103/d4dw-2v6k}}.

\bibitem{Sakai2026jinst}
K.~Sakai, K.~Tomita, T.~Hoshi, R.~Yasuhara, {Monte Carlo} simulation method for
  incoherent {Thomson} scattering spectra from arbitrary electron distribution
  functions, Journal of Instrumentation 21 (2026) C01002.
\newblock \href {https://doi.org/10.1088/1748-0221/21/01/C01002}
  {\path{doi:10.1088/1748-0221/21/01/C01002}}.

\bibitem{Sakai2026ppcf}
K.~Sakai, K.~Tomita, T.~Hoshi, A.~Nakano, M.~Goto, K.~Nagaoka, R.~Yasuhara,
  Conceptual design of {Thomson} scattering system with high wavelength
  resolution in magnetically confined plasmas for electron phase-space
  measurements, Plasma Physics and Controlled Fusion 68 (2026) 015034.
\newblock \href {https://doi.org/10.1088/1361-6587/ae37a4}
  {\path{doi:10.1088/1361-6587/ae37a4}}.

\bibitem{Kobayashi2026}
T.~Kobayashi, T.~Hoshi, A.~Nakano, M.~Yoshinuma, K.~Ida, {Bayesian} modeling of
  charge exchange recombination spectroscopy in fusion plasmas, Review of
  Scientific Instruments 97 (2026) 063506.
\newblock \href {https://doi.org/10.1063/5.0303397}
  {\path{doi:10.1063/5.0303397}}.

\bibitem{Kin2026pfr}
F.~Kin, T.~Kobayashi, T.~Tokuzawa, D.~Nishimura, T.~Hoshi, M.~Sasaki,
  K.~Nagasaki, S.~Inagaki, Reconstruction of electron temperature profiles
  using {Bayesian} analysis on electron cyclotron emission diagnostics in
  {CHD}, Plasma and Fusion Research 21 (2026) 1402023.
\newblock \href {https://doi.org/10.1585/pfr.21.1402023}
  {\path{doi:10.1585/pfr.21.1402023}}.

\bibitem{Kin2026ppcf}
F.~Kin, S.~Kobayashi, A.~Matsuyama, F.~Cai, T.~Kobayashi, T.~Hoshi, S.~Satake,
  S.~Kado, S.~I. Inagaki, S.~Konoshima, T.~Mizuuchi, S.~Inagaki, K.~Nagasaki,
  Development of systematic exploration of magnetic configurations in
  {Heliotron J}, Plasma Physics and Controlled Fusion 68 (2026) 085050.
\newblock \href {https://doi.org/10.1088/1361-6587/ae9beb}
  {\path{doi:10.1088/1361-6587/ae9beb}}.

\bibitem{Kohashi2022}
H.~Kohashi, H.~Iwamoto, T.~Fukaya, Y.~Yamamoto, T.~Hoshi, Performance
  prediction of massively parallel computation by {Bayesian} inference, JSIAM
  Letters 14 (2022) 13--16.
\newblock \href {https://doi.org/10.14495/jsiaml.14.13}
  {\path{doi:10.14495/jsiaml.14.13}}.

\bibitem{Kudo2024}
S.~Kudo, Y.~Yamamoto, T.~Hoshi, A fast and efficient computation method for
  reflective diffraction simulations, Computer Physics Communications 296
  (2024) 109029.
\newblock \href {https://doi.org/10.1016/j.cpc.2023.109029}
  {\path{doi:10.1016/j.cpc.2023.109029}}.

\bibitem{Motoyama2022physbo}
Y.~Motoyama, R.~Tamura, K.~Yoshimi, K.~Terayama, T.~Ueno, K.~Tsuda, Bayesian
  optimization package: {PHYSBO}, Computer Physics Communications 278 (2022)
  108405.
\newblock \href {https://doi.org/10.1016/j.cpc.2022.108405}
  {\path{doi:10.1016/j.cpc.2022.108405}}.

\bibitem{Motoyama2026physbo}
Y.~Motoyama, K.~Yoshimi, T.~Aoyama, K.~Terayama, K.~Tsuda, R.~Tamura, Update of
  physbo: Improving usability and portability of bayesian optimization for
  physics and materials research, Computer Physics Communications 329 (2026)
  110368.
\newblock \href {https://doi.org/10.1016/j.cpc.2026.110368}
  {\path{doi:10.1016/j.cpc.2026.110368}}.

\bibitem{Motoyama2022materiapps}
Y.~Motoyama, K.~Yoshimi, T.~Kato, S.~Todo, {MateriApps LIVE!} and {MateriApps
  Installer}: Environment for starting and scaling up materials science
  simulations, SoftwareX 20 (2022) 101210.
\newblock \href {https://doi.org/10.1016/j.softx.2022.101210}
  {\path{doi:10.1016/j.softx.2022.101210}}.

\bibitem{ForemanMackey2013}
D.~Foreman-Mackey, D.~W. Hogg, D.~Lang, J.~Goodman, emcee: The {MCMC} hammer,
  Publications of the Astronomical Society of the Pacific 125 (2013) 306--312.
\newblock \href {https://doi.org/10.1086/670067} {\path{doi:10.1086/670067}}.

\bibitem{AbrilPla2023}
O.~Abril-Pla, V.~Andreani, C.~Carroll, L.~Dong, C.~J. Fonnesbeck, M.~Kochurov,
  R.~Kumar, J.~Lao, C.~C. Luhmann, O.~A. Martin, M.~Osthege, R.~Vieira,
  T.~Wiecki, R.~Zinkov, {PyMC}: a modern, and comprehensive probabilistic
  programming framework in {Python}, PeerJ Computer Science 9 (2023) e1516.
\newblock \href {https://doi.org/10.7717/peerj-cs.1516}
  {\path{doi:10.7717/peerj-cs.1516}}.

\bibitem{Akiba2019}
T.~Akiba, S.~Sano, T.~Yanase, T.~Ohta, M.~Koyama, Optuna: A next-generation
  hyperparameter optimization framework, in: Proceedings of the 25th ACM SIGKDD
  International Conference on Knowledge Discovery \& Data Mining, 2019, pp.
  2623--2631.
\newblock \href {https://doi.org/10.1145/3292500.3330701}
  {\path{doi:10.1145/3292500.3330701}}.

\bibitem{NelsonPrescott2019}
A.~R.~J. Nelson, S.~W. Prescott, refnx: neutron and {X}-ray reflectometry
  analysis in {Python}, Journal of Applied Crystallography 52 (2019) 193--200.
\newblock \href {https://doi.org/10.1107/S1600576718017296}
  {\path{doi:10.1107/S1600576718017296}}.

\bibitem{GlavicBjorck2022}
A.~Glavic, M.~Bj{\"o}rck, {GenX 3}: the latest generation of an established
  tool, Journal of Applied Crystallography 55 (2022) 1063--1071.
\newblock \href {https://doi.org/10.1107/S1600576722006653}
  {\path{doi:10.1107/S1600576722006653}}.

\bibitem{Virtanen2020}
P.~Virtanen, R.~Gommers, T.~E. Oliphant, et~al., {SciPy} 1.0: fundamental
  algorithms for scientific computing in {Python}, Nature Methods 17 (2020)
  261--272.
\newblock \href {https://doi.org/10.1038/s41592-019-0686-2}
  {\path{doi:10.1038/s41592-019-0686-2}}.

\bibitem{Dalcin2021}
L.~Dalcin, Y.-L.~L. Fang, mpi4py: Status update after 12 years of development,
  Computing in Science \& Engineering 23~(4) (2021) 47--54.
\newblock \href {https://doi.org/10.1109/MCSE.2021.3083216}
  {\path{doi:10.1109/MCSE.2021.3083216}}.

\bibitem{Himmelblau1972}
D.~M. Himmelblau, Applied Nonlinear Programming, McGraw-Hill, 1972.

\bibitem{NelderMead1965}
J.~A. Nelder, R.~Mead, A simplex method for function minimization, The Computer
  Journal 7 (1965) 308--313.
\newblock \href {https://doi.org/10.1093/comjnl/7.4.308}
  {\path{doi:10.1093/comjnl/7.4.308}}.

\bibitem{Sobol1967}
I.~M. Sobol', On the distribution of points in a cube and the approximate
  evaluation of integrals, USSR Computational Mathematics and Mathematical
  Physics 7 (1967) 86--112.
\newblock \href {https://doi.org/10.1016/0041-5553(67)90144-9}
  {\path{doi:10.1016/0041-5553(67)90144-9}}.

\bibitem{Halton1960}
J.~H. Halton, On the efficiency of certain quasi-random sequences of points in
  evaluating multi-dimensional integrals, Numerische Mathematik 2 (1960)
  84--90.
\newblock \href {https://doi.org/10.1007/BF01386213}
  {\path{doi:10.1007/BF01386213}}.

\bibitem{McKay1979}
M.~D. McKay, R.~J. Beckman, W.~J. Conover, A comparison of three methods for
  selecting values of input variables in the analysis of output from a computer
  code, Technometrics 21 (1979) 239--245.
\newblock \href {https://doi.org/10.1080/00401706.1979.10489755}
  {\path{doi:10.1080/00401706.1979.10489755}}.

\bibitem{Owen1998}
A.~B. Owen, Scrambling {Sobol'} and {Niederreiter--Xing} points, Journal of
  Complexity 14 (1998) 466--489.
\newblock \href {https://doi.org/10.1006/jcom.1998.0487}
  {\path{doi:10.1006/jcom.1998.0487}}.

\bibitem{StornPrice1997}
R.~Storn, K.~Price, Differential evolution -- a simple and efficient heuristic
  for global optimization over continuous spaces, Journal of Global
  Optimization 11 (1997) 341--359.
\newblock \href {https://doi.org/10.1023/A:1008202821328}
  {\path{doi:10.1023/A:1008202821328}}.

\bibitem{Endres2018}
S.~C. Endres, C.~Sandrock, W.~W. Focke, A simplicial homology algorithm for
  {Lipschitz} optimisation, Journal of Global Optimization 72 (2018) 181--217.
\newblock \href {https://doi.org/10.1007/s10898-018-0645-y}
  {\path{doi:10.1007/s10898-018-0645-y}}.

\bibitem{Jones1993}
D.~R. Jones, C.~D. Perttunen, B.~E. Stuckman, Lipschitzian optimization without
  the {Lipschitz} constant, Journal of Optimization Theory and Applications 79
  (1993) 157--181.
\newblock \href {https://doi.org/10.1007/BF00941892}
  {\path{doi:10.1007/BF00941892}}.

\bibitem{Xiang1997}
Y.~Xiang, D.~Y. Sun, W.~Fan, X.~G. Gong, Generalized simulated annealing
  algorithm and its application to the {Thomson} model, Physics Letters A 233
  (1997) 216--220.
\newblock \href {https://doi.org/10.1016/S0375-9601(97)00474-X}
  {\path{doi:10.1016/S0375-9601(97)00474-X}}.

\bibitem{Sozykin2022}
K.~Sozykin, A.~Chertkov, R.~Schutski, A.-H. Phan, A.~S. Cichocki, I.~Oseledets,
  {TTOpt}: A maximum volume quantized tensor train-based optimization and its
  application to reinforcement learning, in: Advances in Neural Information
  Processing Systems, Vol.~35, Curran Associates, Inc., 2022, pp. 26052--26065.
\newblock \href {https://doi.org/10.52202/068431-1889}
  {\path{doi:10.52202/068431-1889}}.

\bibitem{Oseledets2011}
I.~V. Oseledets, Tensor-train decomposition, SIAM Journal on Scientific
  Computing 33 (2011) 2295--2317.
\newblock \href {https://doi.org/10.1137/090752286}
  {\path{doi:10.1137/090752286}}.

\bibitem{OseledetsTyrtyshnikov2010}
I.~Oseledets, E.~Tyrtyshnikov, {TT}-cross approximation for multidimensional
  arrays, Linear Algebra and its Applications 432 (2010) 70--88.
\newblock \href {https://doi.org/10.1016/j.laa.2009.07.024}
  {\path{doi:10.1016/j.laa.2009.07.024}}.

\bibitem{Goreinov2010}
S.~A. Goreinov, I.~V. Oseledets, D.~V. Savostyanov, E.~E. Tyrtyshnikov, N.~L.
  Zamarashkin, How to find a good submatrix, in: Matrix Methods: Theory,
  Algorithms and Applications, World Scientific, 2010, pp. 247--256.
\newblock \href {https://doi.org/10.1142/9789812836021_0015}
  {\path{doi:10.1142/9789812836021_0015}}.

\bibitem{WalesDoye1997}
D.~J. Wales, J.~P.~K. Doye, Global optimization by basin-hopping and the lowest
  energy structures of {Lennard-Jones} clusters containing up to 110 atoms, The
  Journal of Physical Chemistry A 101 (1997) 5111--5116.
\newblock \href {https://doi.org/10.1021/jp970984n}
  {\path{doi:10.1021/jp970984n}}.

\bibitem{SwendsenWang1986}
R.~H. Swendsen, J.-S. Wang, Replica {Monte Carlo} simulation of spin-glasses,
  Physical Review Letters 57 (1986) 2607--2609.
\newblock \href {https://doi.org/10.1103/PhysRevLett.57.2607}
  {\path{doi:10.1103/PhysRevLett.57.2607}}.

\bibitem{HukushimaNemoto1996}
K.~Hukushima, K.~Nemoto, Exchange {Monte Carlo} method and application to spin
  glass simulations, Journal of the Physical Society of Japan 65 (1996)
  1604--1608.
\newblock \href {https://doi.org/10.1143/JPSJ.65.1604}
  {\path{doi:10.1143/JPSJ.65.1604}}.

\bibitem{HukushimaIba2003}
K.~Hukushima, Y.~Iba, Population annealing and its application to a spin glass,
  in: AIP Conference Proceedings, Vol. 690, 2003, pp. 200--206.
\newblock \href {https://doi.org/10.1063/1.1632130}
  {\path{doi:10.1063/1.1632130}}.

\bibitem{Machta2010}
J.~Machta, Population annealing with weighted averages: A {Monte Carlo} method
  for rough free-energy landscapes, Physical Review E 82 (2010) 026704.
\newblock \href {https://doi.org/10.1103/PhysRevE.82.026704}
  {\path{doi:10.1103/PhysRevE.82.026704}}.

\bibitem{Neal2001}
R.~M. Neal, Annealed importance sampling, Statistics and Computing 11 (2001)
  125--139.
\newblock \href {https://doi.org/10.1023/A:1008923215028}
  {\path{doi:10.1023/A:1008923215028}}.

\bibitem{Weigel2021}
M.~Weigel, L.~Barash, L.~Shchur, W.~Janke, Understanding population annealing
  {Monte Carlo} simulations, Physical Review E 103 (2021) 053301.
\newblock \href {https://doi.org/10.1103/PhysRevE.103.053301}
  {\path{doi:10.1103/PhysRevE.103.053301}}.

\bibitem{Hanada2022}
T.~Hanada, Y.~Motoyama, K.~Yoshimi, T.~Hoshi, {sim-trhepd-rheed} -- open-source
  simulator of total-reflection high-energy positron diffraction ({TRHEPD}) and
  reflection high-energy electron diffraction ({RHEED}), Computer Physics
  Communications 277 (2022) 108371.
\newblock \href {https://doi.org/10.1016/j.cpc.2022.108371}
  {\path{doi:10.1016/j.cpc.2022.108371}}.

\bibitem{Voegeli2006}
W.~Voegeli, K.~Akimoto, T.~Aoyama, K.~Sumitani, S.~Nakatani, H.~Tajiri,
  T.~Takahashi, Y.~Hisada, S.~Mukainakano, X.~Zhang, H.~Sugiyama, H.~Kawata,
  Structure of the {SiC(0001)} $3\times3$ reconstruction studied by surface
  x-ray diffraction, Applied Surface Science 252 (2006) 5259--5262.
\newblock \href {https://doi.org/10.1016/j.apsusc.2005.12.019}
  {\path{doi:10.1016/j.apsusc.2005.12.019}}.

\bibitem{Ankudinov1998}
A.~L. Ankudinov, B.~Ravel, J.~J. Rehr, S.~D. Conradson, Real-space
  multiple-scattering calculation and interpretation of x-ray-absorption
  near-edge structure, Physical Review B 58~(12) (1998) 7565--7576.
\newblock \href {https://doi.org/10.1103/PhysRevB.58.7565}
  {\path{doi:10.1103/PhysRevB.58.7565}}.

\bibitem{Artrith2016}
N.~Artrith, A.~Urban, An implementation of artificial neural-network potentials
  for atomistic materials simulations: Performance for {TiO2}, Computational
  Materials Science 114 (2016) 135--150.
\newblock \href {https://doi.org/10.1016/j.commatsci.2015.11.047}
  {\path{doi:10.1016/j.commatsci.2015.11.047}}.

\bibitem{Akai1989}
H.~Akai, Fast {Korringa-Kohn-Rostoker} coherent potential approximation and its
  application to {FCC} {Ni-Fe} systems, Journal of Physics: Condensed Matter
  1~(43) (1989) 8045.
\newblock \href {https://doi.org/10.1088/0953-8984/1/43/006}
  {\path{doi:10.1088/0953-8984/1/43/006}}.

\bibitem{Yoshimi2026data}
T.~Aoyama, K.~Yoshimi, Y.~Motoyama, A.~Nakano, H.~Li, Y.~Konishi, T.~Hoshi,
  \href{https://isspns-gitlab.issp.u-tokyo.ac.jp/odatse-dev/odatse-v4-paper-data}{Snapshot
  data for ``{ODAT-SE}: A modular and parallel framework for simulation-based
  inverse problems''}, ISSP Data Repository, The University of Tokyo (2026).
\newline\urlprefix\url{https://isspns-gitlab.issp.u-tokyo.ac.jp/odatse-dev/odatse-v4-paper-data}

\bibitem{JamilYang2013}
M.~Jamil, X.-S. Yang, A literature survey of benchmark functions for global
  optimisation problems, International Journal of Mathematical Modelling and
  Numerical Optimisation 4~(2) (2013) 150--194.
\newblock \href {https://doi.org/10.1504/IJMMNO.2013.055204}
  {\path{doi:10.1504/IJMMNO.2013.055204}}.

\bibitem{SurjanovicBingham2013}
S.~Surjanovic, D.~Bingham, Virtual library of simulation experiments: Test
  functions and datasets, \url{http://www.sfu.ca/~ssurjano}, retrieved
  September 23, 2026 (2013).

\bibitem{Ichimiya1987}
A.~Ichimiya, {RHEED} intensity analysis of {Si(111)7$\times$7} at one-beam
  condition, Surface Science 192 (1987) L893--L898.
\newblock \href {https://doi.org/10.1016/S0039-6028(87)81122-6}
  {\path{doi:10.1016/S0039-6028(87)81122-6}}.

\bibitem{Shirasawa2006}
T.~Shirasawa, S.~Mizuno, H.~Tochihara, Structural analysis of the
  {c(4$\times$2)} reconstruction in {Si(001)} and {Ge(001)} surfaces by
  low-energy electron diffraction, Surface Science 600 (2006) 815--819.
\newblock \href {https://doi.org/10.1016/j.susc.2005.11.031}
  {\path{doi:10.1016/j.susc.2005.11.031}}.

\bibitem{Yoshimi2025}
K.~Yoshimi, Y.~Motoyama, T.~Aoyama, M.~Kawamura, N.~Kawashima, Project for
  advancement of software usability in materials science, Science and
  Technology of Advanced Materials: Methods 5~(1) (2025) 2564055.
\newblock \href {https://doi.org/10.1080/27660400.2025.2564055}
  {\path{doi:10.1080/27660400.2025.2564055}}.

\bibitem{RasmussenWilliams2005}
C.~E. Rasmussen, C.~K.~I. Williams, Gaussian Processes for Machine Learning,
  The MIT Press, 2006.
\newblock \href {https://doi.org/10.7551/mitpress/3206.001.0001}
  {\path{doi:10.7551/mitpress/3206.001.0001}}.

\bibitem{Pedregosa2011}
F.~Pedregosa, G.~Varoquaux, A.~Gramfort, V.~Michel, B.~Thirion, O.~Grisel,
  M.~Blondel, P.~Prettenhofer, R.~Weiss, V.~Dubourg, J.~Vanderplas, A.~Passos,
  D.~Cournapeau, M.~Brucher, M.~Perrot, {\'E}.~Duchesnay, Scikit-learn: Machine
  learning in {Python}, Journal of Machine Learning Research 12 (2011)
  2825--2830.

\bibitem{Lucas2014}
A.~Lucas, Ising formulations of many {NP} problems, Frontiers in Physics 2
  (2014) 5.
\newblock \href {https://doi.org/10.3389/fphy.2014.00005}
  {\path{doi:10.3389/fphy.2014.00005}}.

\bibitem{Bussieck2003}
M.~R. Bussieck, A.~S. Drud, A.~Meeraus, {MINLPLib} -- a collection of test
  models for mixed-integer nonlinear programming, INFORMS Journal on Computing
  15 (2003) 114--119.
\newblock \href {https://doi.org/10.1287/ijoc.15.1.114.15159}
  {\path{doi:10.1287/ijoc.15.1.114.15159}}.

\end{thebibliography}
\end{document}